\documentclass[onefignum,onetabnum]{siamart171218} 

\usepackage[utf8]{inputenc} 

\usepackage{mathtools}
\usepackage{color}
\usepackage{amsfonts}
\usepackage{graphicx}
\usepackage{epstopdf}
\usepackage{algorithmic}
\usepackage{subcaption}
\usepackage{tabularx}
\usepackage{tikz,pgfplots}
\newlength\figheight
\newlength\figwidth

\ifpdf
  \DeclareGraphicsExtensions{.eps,.pdf,.png,.jpg}
\else
  \DeclareGraphicsExtensions{.eps}
\fi

\newsiamremark{remark}{Remark}
\newsiamremark{hypothesis}{Hypothesis}
\crefname{hypothesis}{Hypothesis}{Hypotheses}
\newsiamthm{claim}{Claim}

\headers{Gradient-based dimension reduction}{O.~Zahm, P.~Constantine, C.~Prieur, Y.~Marzouk}

\title{Gradient-based dimension reduction of multivariate vector-valued functions
\thanks{Submitted to the editors on 20 October  2018.
}}

\author{
  Olivier Zahm\thanks{Univ. Grenoble Alpes, Inria, CNRS, Grenoble INP (Institute of Engineering Univ. Grenoble Alpes), LJK, 38000 Grenoble, France \email{olivier.zahm@inria.fr}}
  \and
  Paul Constantine\thanks{University of Colorado Boulder, \email{paul.constantine@colorado.edu}}
  \and
  Cl\'ementine Prieur\thanks{Univ. Grenoble Alpes, CNRS, Inria, Grenoble INP (Institute of Engineering Univ. Grenoble Alpes), LJK, Grenoble 38000, France \email{clementine.prieur@univ-grenoble-alpes.fr}}
  \and
  Youssef Marzouk\thanks{Massachusetts Institute of Technology, \email{ymarz@mit.edu}}
}

\usepackage{amsopn}

\DeclareMathOperator{\trace}{trace}
\DeclareMathOperator{\Var}{Var}

\ifpdf
\hypersetup{
  pdftitle={Gradient-based dimension reduction of multivariate vector-valued functions},
  pdfauthor={O.~Zahm, P.~Constantine, C.~Prieur, Y.~Marzouk}
}
\fi

\begin{document}

\maketitle

\begin{abstract}
  Multivariate functions encountered in high-dimensional uncertainty quantification problems often vary 
 along a few dominant directions in the input parameter space. We propose a gradient-based method for detecting these directions and using them to construct {ridge approximations} of such functions, in a setting where the functions are vector-valued (e.g., taking values in $\mathbb{R}^n$).
The methodology consists of minimizing an upper bound on the approximation error, obtained by {subspace Poincaré inequalities}.
We provide a thorough mathematical analysis in the case where the parameter space is equipped with a Gaussian probability measure.
The resulting method generalizes the notion of active subspaces associated with scalar-valued functions. 
A numerical illustration shows that using gradients of the function yields effective dimension reduction.
We also show how the choice of norm on the codomain of the function has an impact on the function's low-dimensional approximation.

\end{abstract}

\begin{keywords}
  High-dimensional function approximation, dimension reduction, active subspace, ridge approximation, Karhunen-Loève decomposition, Poincaré inequality, Sobol' indices.
\end{keywords}

\begin{AMS}
  41A30, 41A63, 65D15
\end{AMS}

\section{Introduction}

Many problems that arise in uncertainty quantification---e.g., integrating or approximating multivariate functions---suffer from the curse of dimensionality: the complexity of algorithms grows dramatically (typically exponentially) with the dimension of the input parameter space. One approach to alleviate this curse is to identify and exploit some notion of low-dimensional structure. For example, the function of interest might vary primarily along a few directions of the input parameter space while being (almost) constant in the other directions. In this case, we say that the problem has a \emph{low intrinsic dimension}; algorithms for quantifying uncertainty can then focus on these important directions to reduce the overall cost.

%
A common and simple approach for parameter space dimension reduction is the {truncated Karhunen-Lo\`{e}ve decomposition}~\cite{Schwab2006}, closely related to {principal component analysis}~\cite{Jolliffe2002}. These techniques exploit the correlation structure of the function's input space (specifically, decay in the spectrum of the covariance of the input measure). However, more effective dimension reduction is possible with techniques that exploit not only input correlations but also the structure of the input-output map itself.
%
One way to reduce the input space dimension is to determine the non-influential input parameters (or factors) and to fix them to some arbitrary value. Factor fixing (see, e.g., \cite{saltar04}) is often a goal of \emph{global sensitivity analysis}~\cite{Saltelli2008,iooss2015review}. 
For independent inputs, total Sobol' indices \cite{saltar04} are a popular way to address the factor fixing problem, as they measure the total impact that each variable (or each group of variables) has on the variance of the output.
Estimating these indices can be computationally challenging, however; see for instance \cite{sobol1993sensitivity,homma1996importance,sobol2001global}. 
Alternative screening procedures based on \emph{derivative-based global sensitivity measures} (DGSM) have been proposed in \cite{kucherenko2009derivative,kucherenko2009monte}.
These indices are defined as integrals of squared derivatives of the model output.
If the numerical implementation of a model permits easy computation of the derivatives (for instance using the \emph{adjoint method}, see \cite{plessix2006review}), these indices can be estimated with reasonable computational cost.
There are interesting links between DGSM and Sobol' indices.
For instance, assuming the inputs are independent, one can bound the total Sobol' indices by the DGSM up to some Poincaré constant that depends on the probability distribution of the parameters (see, e.g.,  \cite{kucherenko2009derivative,lamboni2013derivative} and \cite{roustant2017poincare} for a recent detailed analysis). 
Yet the factor fixing setting is somewhat restrictive, in that functions often vary most prominently in directions that are \emph{not aligned} with the coordinate axes corresponding to the original inputs.

Closely related to derivative-based screening are \emph{active subspaces}, described in \cite{russi2010uncertainty,Constantine2015,Constantine2014a}. 
Active subspaces are defined as the leading eigenspaces of the second moment matrix of the function's gradient, the diagonal of which contains the DGSM. These eigenspaces are not necessarily aligned with the canonical coordinates, and hence are able to identify linear combinations of the input parameters along which the function varies the most. In this sense, they generalize coordinate-aligned derivative-based global sensitivity analysis.  Active subspaces have been used in a wide range of science and engineering models~\cite{Lukaczyk2014,Constantine2017a, Jefferson2017}.  Connections between Sobol’ indices, DGSM, and active subspaces for scalar-valued functions are explored in \cite{constantine2017global}.

Global sensitivity analysis and active subspaces have primarily been focused on scalar-valued functions, as in models with a single output quantity of interest. In the presence of multiple outputs of interest, as is the case in many practical applications, new approaches are needed. Aggregated Sobol' indices for multiple outputs or functional outputs have been introduced in \cite{lamboni2011multivariate}, and further studied in \cite{gamboa:hal-00800847,gamboa:hal-00881112}. In the context of active subspaces, one could try to identify important input parameter directions for each output and then combine all those directions, as in \cite{ji2018shared}. But it is not clear how to interpret or even best perform such a combination step.

\subsection{Contribution}
In this paper, we propose a methodology for detecting and exploiting the low intrinsic dimension that a given multivariate function might have. We formulate our approach as a \emph{controlled approximation problem}, seeking a certified upper bound for the error in a \emph{ridge} approximation of the original function. With this approximation perspective, our methodology extends naturally to the case of vector-valued functions---for instance, functions with multiple real-valued outputs. Specifically, given a function of interest 
\begin{equation*}
 x\mapsto f(x_1,\hdots,x_d) \in V ,
\end{equation*}
where $V$ is a vector space, the problem is to find an approximation of $f$ by a function of fewer variables, say $y\mapsto g(y_1,\hdots,y_r)$ with $r\ll d$ where $y=h(x)$ depends linearly on $x$.
Thus, given a user-defined tolerance $\varepsilon$, we seek a \emph{linear} function $h$ such that
\begin{equation}\label{eq:INTRO}
 \| f-g\circ h \| \leq \varepsilon
 \qquad\text{where}\quad\left\{
 \begin{array}{l}
  \mathbb{R}^d \xrightarrow{\quad\quad f \quad\quad} V \\
  \mathbb{R}^{d} \overset{h}{\longrightarrow} {\mathbb{R}^{r}} \overset{g}{\longrightarrow} V ,
 \end{array}\right
 .
\end{equation}
holds for some function $g$, where $\|\cdot\|$ is a norm chosen depending on the application.
Approximations of the form of $g\circ h$ are called \emph{ridge functions} \cite{Pinkus2015}.
If such an approximation exists with $r\ll d$, we say that $f$ has a \emph{low effective dimension} $r=r(\varepsilon)$, and $(y_1,\hdots,y_r)=h(x)$ correspond to the \emph{active} (or \emph{explanatory}) variables.
To solve this controlled approximation problem, we use Poincaré-type inequalities to derive an upper bound on the error. This bound, defined by means of gradients (or Jacobians) of $f$, admits a simple expression and can be analytically minimized with respect to $h$ and $g$ for any fixed $r$.
By choosing $r$ such that the minimized error bound is below the prescribed tolerance $\varepsilon$, we obtain an approximation of $f$ whose error is controlled. We also show that, for scalar-valued functions $f$, the minimizer of the bound corresponds to the active subspace approach proposed in \cite{Constantine2014a}.

In our analysis we assume that the parameter domain is equipped with a Gaussian probability measure, and we define the norm $\|\cdot\|$ in \eqref{eq:INTRO} as the corresponding weighted norm. Thus \eqref{eq:INTRO} becomes an approximation problem for $f$ in the \emph{mean-squared} sense. The Gaussian measure need not be standard: it can have non-zero mean and non-identity covariance matrix.  By allowing the latter, we will show that the notion of a low effective dimension also depends on the input covariance matrix itself. Furthermore, having non-standard Gaussian measures enables us to compare our approach with the truncated Karhunen-Lo\`{e}ve decomposition, which also exploits the spectral properties of the parameter covariance matrix. 
The Gaussian assumption primarily permits us to simplify our analysis. One can consider other probability measures as long as they satisfy the so-called \emph{subspace Poincaré inequality} described later in the paper, which is the key argument of our method. Explicit generalizations of  this inequality to non-Gaussian measures are given in \cite{zahm2018certified}.

It is important to mention that in actual practice, minimizing the error itself is a much more difficult problem than minimizing the error bound. This is why the proposed strategy is appealing, provided that gradient information from $f$ is available. However, there is no guarantee that the minimizer of the bound is close to the minimizer of the true error.
To illustrate the potential and the limitations of the proposed method, we present (in Section~\ref{sec:analytical}) examples of functions $f$ for which minimizing the bound gives either the minimizer of the error (ideal case) or the maximizer of the error (worst case). In both cases the proposed method still permits us to control the approximation error: it simply does so more efficiently in the first case than in the second.
We also demonstrate our method on a parameterized partial differential equation (see Section~\ref{sec:pde}). This example shows that the resulting ridge approximation depends not only on $f$ but also on the choice of norm on the output space $V$, which in turn defines the function-space norm $\|\cdot\|$  in \eqref{eq:INTRO}.

\medskip

Ridge functions and their approximation properties were extensively studied in the 1980s because of their connection to both projection pursuit regression~\cite{Friedman1981,Diaconis1984,Huber1985} and early neural networks~\cite{Haykin1999}. Recent work has exploited compressed sensing to recover a ridge function from point queries~\cite{Fornasier2012,Cohen2012b}. 
The ridge recovery problem corresponds to the proposed problem setup \eqref{eq:INTRO} with $\varepsilon=0$: the goal is to recover $g$, $h$, and $r$ assuming that $f$ is exactly a ridge function $f=g\circ h$.
In contrast, we do not aim for an exact recovery of $f$, but rather approximation of $f$ by a ridge function up to a prescribed precision $\varepsilon>0$.
Similar recovery problems arise in the statistical regression literature under the name \emph{sufficient dimension reduction}~\cite{Adragni09,Cook1998}.
In this context, the goal is to identify linear combinations in the input space that are \emph{statistically sufficient} to explain the regression response. Among the numerous sufficient dimension reduction techniques that have been proposed, we mention sliced inverse regression \cite{li1991sliced}, sliced average variance estimation \cite{cookweis1991}, and principal Hessian directions \cite{li1992principal}. 
In \cite{samarov1993exploring}, gradient information is used to explore the underlying regression structure by means of \emph{average derivative functionals}, estimated nonparametrically via kernels.  
Concerning dimension reduction in regression with vector-valued responses, a broad literature has also emerged more recently. We refer to \cite{zhu2010dimension} and references therein (see also \cite{li2003dimension,saracco2005asymptotics,barreda2007some}).
Broadly, and in contrast with the approach proposed here, these regression analyses are concerned with estimation from a given data set, and thus rely on statistical assessments of the error. 
\medskip

The rest of this paper is organized as follows. Section \ref{sec2} describes our dimension reduction methodology, deriving an upper bound on the error and an explicit construction for its minimizer, yielding a controlled ridge approximation of a vector-valued function. 
Section~\ref{sec:KL} compares the proposed method with the truncated Karhunen-Lo\`{e}ve decomposition, and Section \ref{sec:DGSM} discusses its relationship with sensitivity analysis. In Section \ref{sec4} we demonstrate our method on various analytical and numerical examples. 
Proofs of the main results are deferred to Appendix \ref{sec:proofs}.

\section{Dimension reduction of the input parameter space} \label{sec2}
Throughout the paper, the algebraic space $\mathbb{R}^d$ refers to a parameter space of dimension $d\gg1$. The Borel sets of $\mathbb{R}^d$ are denoted by $\mathcal{B}(\mathbb{R}^d)$ and we let $\mu=\mathcal{N}(m,\Sigma)$ be the Gaussian probability measure on $\mathbb{R}^d$ with mean $m\in\mathbb{R}^d$ and covariance $\Sigma\in\mathbb{R}^{d\times d}$, which is assumed to be non-singular. We let $V=\mathbb{R}^n$ be an algebraic space endowed with a norm $\|\cdot\|_V$ associated with a scalar product $(\cdot,\cdot)_V$ defined by $(v,w)_V = v^T R_V w$ for any $v,w\in V$, where  $R_V\in\mathbb{R}^{n\times n}$ is a symmetric positive definite matrix.
We denote by
$$
 \mathcal{H} = L^2( \mathbb{R}^d , \mathcal{B}(\mathbb{R}^d) ,\mu ; V )  \, ,
$$
the Hilbert space which contains all the measurable functions $v:\mathbb{R}^d\rightarrow V$ such that $\| v \|_\mathcal{H} <\infty$, where $\| \cdot \|_\mathcal{H} $ is the norm associated with the scalar product $(\cdot,\cdot)_\mathcal{H}$ defined by
$$
 (u,v)_\mathcal{H} = \int (u(x),v(x))_V \, \mathrm{d}\mu(x)  \, ,
$$
for any $u,v\in\mathcal{H}$.

Ridge functions are functions of the form $g\circ h$ where $h : \mathbb{R}^d \rightarrow \mathbb{R}^r$ is a linear function and where $g : \mathbb{R}^r \rightarrow V$ is a measurable function, sometimes called the \emph{profile} of the ridge function; see \cite{Mayer2015}.
Ridge functions are essentially functions that are constant along a subspace (the kernel of $h$).
In this paper we will use the following parametrization of ridge functions,
\begin{equation}\label{eq:RigdeFunction}
 x \mapsto g(P_rx)  ,
\end{equation}
where $P_r \in\mathbb{R}^{d\times d}$ is a rank-$r$ projector and $g : \mathbb{R}^d \rightarrow V$ is a measurable function.
Notice that $g(P_rx) = g(P_ry)$ whenever $x-y\in\text{Ker}(P_r)$, which means that the function \eqref{eq:RigdeFunction} is constant along the kernel of the projector, and thus is a ridge function.\footnote{One can easily show that any function as in \eqref{eq:RigdeFunction} can be written as $g' \circ h$ for some linear $h:\mathbb{R}^d\rightarrow\mathbb{R}^r$ and some measurable $g':\mathbb{R}^r\rightarrow V$, and vice versa.}

We consider the problem of finding a controlled approximation of a function $f \in \mathcal{H}$ by a ridge function.
Given a prescribed tolerance $\varepsilon \geq 0$, the problem consists in finding $g$ and $P_r$ such that
\begin{equation}\label{eq:ControlledApprox}
 \| f- g\circ P_r \|_{\mathcal{H}} \leq \varepsilon \, .
\end{equation}
The choice $P_r=I_d$ (the identity matrix) and $g=f$ in \eqref{eq:RigdeFunction} yields a trivial solution. But in that case, the rank of $P_r$ is equal to $d$ and there is no dimension reduction. Thus, in order to make this problem meaningful, we want $r=\text{rank}(P_r)$ to be less than $d$, ideally $r\ll d$. 

\begin{remark}
An equivalent formulation of the problem is the following. Given a tolerance $\varepsilon>0$, we want to find a Borel function $g:\mathbb{R}^d \rightarrow V $ and a low-rank projector $P_r\in\mathbb{R}^{d\times d}$ such that
 $$
  \mathbb{E}\big(\| f(X)-g(P_rX) \|_V^2\big) \leq \varepsilon^2  ,
 $$
where $X\sim\mathcal{N}(m,\Sigma)$ is a random vector and where $\mathbb{E}(\cdot)$ denotes the mathematical expectation. If $\varepsilon^2 \ll \Var(f(X))=\mathbb{E}(\|f(X) - \mathbb{E}(f(X))\|_V^2)$, the statistical interpretation is that the random variable $X_r = P_rX$ is an explanatory variable for $f(X)$, in the sense that most of the variance of $f(X)$ can be explained by $X_r$.
 
\end{remark}

\subsection{Optimal profile for the ridge function}\label{sec:OptimalG}

In this section, we assume that the projector $P_r$ is given. We denote by
$$ 
 \mathcal{H}_{P_r} = L^2( \mathbb{R}^d , \sigma(P_r) ,\mu ; V ) ,
$$
the space containing all the $\sigma(P_r)$-measurable functions $v:\mathbb{R}^d\rightarrow V$ such that $\|v\|_{\mathcal{H}}<\infty$. Here $\sigma(P_r)$ is the $\sigma$-algebra generated by $P_r$. By the Doob--Dynkin lemma, see for example Lemma 1.13 in \cite{Mikosch1998}, the set of all $\sigma(P_r)$-measurable functions is exactly the set of the functions of the form $x\mapsto g(P_rx)$ for some Borel function $g$, so that 
\begin{equation}\label{eq:Doob}
 \mathcal{H}_{P_r} = \{ g\circ P_r ~|~ g:\mathbb{R}^d\rightarrow V \text{, Borel function} \} \cap \mathcal{H} .
\end{equation}
Note that $\mathcal{H}_{P_r}$ is a closed subspace in $\mathcal{H}$. Then, for any $f\in\mathcal{H}$, there exists a unique minimizer of $f_r \mapsto \| f - f_r \|_\mathcal{H}$ over $\mathcal{H}_{P_r}$. This minimizer corresponds to the orthogonal projection of $f\in\mathcal{H}$ onto $\mathcal{H}_{P_r}$ and is denoted by $\mathbb{E}_\mu(f|\sigma(P_r) )$. We can write
$$
 \| f- \mathbb{E}_\mu(f|\sigma(P_r) ) \|_{\mathcal{H}}
 =
 \min_{f_r\in \mathcal{H}_{P_r}} \| f- f_r \|_{\mathcal{H}}
 = 
 \min_{\substack{g:\mathbb{R}^d \rightarrow V \\ \text{Borel function}}} \| f-g\circ P_r  \|_{\mathcal{H}} ,
$$
which means that $\mathbb{E}_\mu(f|\sigma(P_r) )$ yields an optimal profile $g$.
Note that $\mathbb{E}_\mu(f|\sigma(P_r) )\in\mathcal{H}_{P_r}$ can be uniquely characterized by the variational equation
\begin{equation}\label{eq:CondExp_VarForm}
\int ( \mathbb{E}_\mu(f|\sigma(P_r)) , h )_V ~\mathrm{d}\mu
= \int ( f , h )_V~\mathrm{d}\mu  \, ,
\end{equation}
for all $h\in\mathcal{H}_{P_r}$. In other words, $\mathbb{E}_\mu(f|\sigma(P_r))$ corresponds to the conditional expectation of $f$ under the distribution $\mu$ given the $\sigma$-algebra $\sigma(P_r)$, which explains the choice of notation. 
The following proposition gives an interesting property on the space $\mathcal{H}_{P_r}$. The proof is given in Appendix \ref{proof:prop:ImportanceKernel}.

\begin{proposition}\label{prop:ImportanceKernel}
 Let $P_r$ and $Q_r$ be two projectors such that $\text{Ker}(P_r) = \text{Ker}(Q_r)$. Then we have
 $
  \mathcal{H}_{P_r} = \mathcal{H}_{Q_r}.
 $
\end{proposition}

Let us recall that a projector is uniquely characterized by both its kernel and its image.\footnote{Of course an \textit{orthogonal} projector (orthogonal with respect to any scalar product) is uniquely characterized either by its kernel or by its image, since the other subspace can be uniquely defined as the orthogonal complement.}
Proposition \ref{prop:ImportanceKernel} shows that $\mathcal{H}_{P_r}$ is invariant with respect to the image of $P_r$, and so is the conditional expectation $\mathbb{E}_\mu(f|\sigma(P_r))$. In particular, the error $P_r\mapsto\| f-\mathbb{E}_\mu(f|\sigma(P_r)) \|_{\mathcal{H}}$ depends only on the kernel of $P_r$. 
This means that, with regard to the initial dimension reduction problem \eqref{eq:ControlledApprox}, the goal is now to find a subspace where the function $f$ does \emph{not} vary.

By Proposition \ref{prop:ImportanceKernel} and without loss of generality, we can assume that $P_r$ is an \textit{orthogonal} projector with respect to an arbitrary scalar product on $\mathbb{R}^d$. In the present context, the natural scalar product to use is the one induced by the precision matrix $\Sigma^{-1}$ of $\mu$, which is $\langle x,y\rangle=x^T\Sigma^{-1}y$ for any $x,y\in\mathbb{R}^d$. The associated norm $\|\cdot\|_{\Sigma^{-1}}$ is such that $\| x \|_{\Sigma^{-1}}^2=x^T\Sigma^{-1} x$ for any $x\in\mathbb{R}^d$.
The projector $P_r$ is $\Sigma^{-1}$-orthogonal if $\langle P_r x , (I_d-P_r)x\rangle=0$ for all $x\in\mathbb{R}^d$, which is equivalent to 
\begin{equation}\label{eq:orthogonalityOfPr}
 \|x\|_{\Sigma^{-1}}^2=\|P_rx\|_{\Sigma^{-1}}^2+\|(I_d-P_r)x\|_{\Sigma^{-1}}^2,
\end{equation}
for all $x\in\mathbb{R}^d$. The following proposition gives a simple expression for the conditional expectation $\mathbb{E}_\mu(f|\sigma(P_r))$, provided $P_r$ satisfies \eqref{eq:orthogonalityOfPr}. 
The proof is given in Appendix \ref{proof:prop:ExplicitCondExp}.

\begin{proposition}\label{prop:ExplicitCondExp}
 Let $\mu=\mathcal{N}(m,\Sigma)$ where $\Sigma\in\mathbb{R}^{d\times d}$ is a non-singular covariance matrix and $f \in \mathcal{H}$. Then for any $\Sigma^{-1}$-orthogonal projector $P_r$ we have
 $$
  \mathbb{E}_\mu(f|\sigma(P_r)) : x\mapsto \mathbb{E}(f(P_r x + (I_d-P_r)Y)) ,
 $$
 where the expectation is taken over the random vector $Y\sim\mu$.
 
\end{proposition}

\subsection{Poincaré-based upper bound for the error}

In this section we show how Poincaré-type inequalities can be used to derive an upper bound for the error. This upper bound holds for any projector and is quadratic in $P_r$ so that it can easily be minimized.
\\

It is well known that the standard Gaussian distribution $\gamma=\mathcal{N}(0,I_d)$ satisfies the \emph{Poincaré inequality}
\begin{equation}\label{eq:Poincare_canonical}
 \int ( h - \mathbb{E}_\gamma(h) )^2 \, \mathrm{d}\gamma
 \leq \int \| \nabla h \|_{2}^2 \, \mathrm{d}\gamma,
\end{equation}
for any continuously differentiable function $h: \mathbb{R}^d\rightarrow \mathbb{R}$, where $\nabla h$ denotes the gradient of $h$ (see for example Theorem 3.20 in \cite{boucheron2013concentration}). Here $\mathbb{E}_\gamma(h) = \int h \, \mathrm{d}\gamma$ and $\|\cdot\|_2 = \sqrt{(\cdot)^T(\cdot)}$ denotes the canonical norm of $\mathbb{R}^d$. As noticed in \cite{Chen1982}, non-standard Gaussian distributions also satisfy a Poincaré inequality. By replacing $h$ by $x\mapsto h( \Sigma^{1/2} x + m )$ in \eqref{eq:Poincare_canonical}, where $\Sigma^{1/2}$ is a symmetric square root of $\Sigma$, we have that $\mu=\mathcal{N}(m,\Sigma)$ satisfies
\begin{equation}\label{eq:Poincare}
 \int ( h - \mathbb{E}_\mu(h) )^2 \, \mathrm{d}\mu
 \leq \int \| \nabla h \|_{\Sigma}^2 \, \mathrm{d}\mu,
\end{equation}
for any continuously differentiable function $h: \mathbb{R}^d\rightarrow \mathbb{R}$, where $\|\cdot\|_{\Sigma}$ is the norm on $\mathbb{R}^d$ such that $\|x\|_{\Sigma}^{2}=x^T\Sigma x$ for all $x\in\mathbb{R}^d$.
The next proposition shows that $\mu$ satisfies another Poincaré-type inequality which we call the \textit{subspace Poincaré inequality}. 
The proof is given in Appendix \ref{proof:prop:subspin}.

\begin{proposition}\label{prop:subspin}
 The probability distribution $\mu=\mathcal{N}(m,\Sigma)$ satisfies 
 \begin{equation}\label{eq:SubspacePoincare}
 \int ( h - \mathbb{E}_\mu(h|\sigma(P_r)) )^2 \, \mathrm{d}\mu
 \leq \int \| (I_d-P_r^T)\nabla h \|_{\Sigma}^2 \, \mathrm{d}\mu ,
\end{equation}
 for any continuously differentiable function $h: \mathbb{R}^d\rightarrow \mathbb{R}$ and for any projector $P_r$.
\end{proposition}

The subspace Poincaré inequality stated in Proposition \ref{prop:subspin} allows us to derive an upper bound for the error $\| f - \mathbb{E}_\mu(f|\sigma(P_r)) \|_{\mathcal{H}}$, as shown by the following proposition, whose proof is given in Appendix \ref{proof:prop:BoundVVfunction}.

\begin{proposition}\label{prop:BoundVVfunction}
 Let $\mu=\mathcal{N}(m,\Sigma)$, where $\Sigma\in\mathbb{R}^{d\times d}$ is a non-singular covariance matrix, and let $f\in\mathcal{H} = L^2( \mathbb{R}^d , \mathcal{B}(\mathbb{R}^d) ,\mu ; V ) $, where $V=\mathbb{R}^n$ is endowed with a norm $\|\cdot\|_V$ such that $\|v\|_V^2=v^TR_Vv$ for some symmetric positive definite matrix $R_V\in\mathbb{R}^{n\times n}$. Furthermore, assume that $f$ is continuously differentiable. Then for any projector $P_r\in\mathbb{R}^{d\times d}$ we have
 \begin{equation}\label{eq:BoundVVfunction}
  \| f - \mathbb{E}_\mu(f|\sigma(P_r)) \|_{\mathcal{H}}^2 \leq {\trace\big( \Sigma(I_d-P_r^T) H (I_d-P_r) \big) },
 \end{equation}
 where $H\in\mathbb{R}^{d\times d}$ is the matrix defined by
 \begin{equation}\label{eq:defH}
  H = \int_{\mathbb{R}^d} (\nabla f(x))^T R_V (\nabla f(x)) ~ \mathrm{d}\mu(x).
 \end{equation}
 Here, $\nabla f(x)\in \mathbb{R}^{n \times d}$ denotes the Jacobian matrix of $f(x)=(f_1(x),\hdots,f_n(x))$ at point $x$ given by
 \begin{equation}\label{eq:defNablaF}
   \nabla f(x) = 
   \left(\begin{array}{ccc}
    \frac{\partial f_1}{\partial x_1} (x)& \cdots & \frac{\partial f_1}{\partial x_d}(x)\\
    \vdots & \ddots & \vdots \\
    \frac{\partial f_n}{\partial x_1} (x)& \cdots & \frac{\partial f_n}{\partial x_d}(x)
   \end{array}\right)
   .
 \end{equation}
\end{proposition}

Note that the matrix $H$ defined in \eqref{eq:defH} depends not only on $f$ but also on the norm $\|\cdot\|_V$ of the output space $V$ via the matrix $R_V$.

\subsection{Minimizing the upper bound}

The following proposition enables minimization of the upper bound in Proposition \ref{prop:BoundVVfunction}.
The proof is given in Appendix \ref{proof:prop:OptimalPr}.

\begin{proposition}\label{prop:OptimalPr}
 Let $\Sigma\in\mathbb{R}^{d\times d}$ be a symmetric positive-definite matrix and $H\in\mathbb{R}^{d\times d}$ a symmetric positive-semidefinite matrix. Denote by $(\lambda_i,v_i)\in\mathbb{R}_{\geq0}\times\mathbb{R}^d$ the $i$-th generalized eigenpair of the matrix pair $(H,\Sigma^{-1})$, meaning $Hv_i = \lambda_i \Sigma^{-1} v_i$ with $\|v_i\|_{\Sigma^{-1}}=1$.
 For any $r \leq d$ we have
 \begin{equation}\label{eq:OptimalPr_residual}
  \min_{\substack{ P_r\in\mathbb{R}^{d\times d} \\ \text{rank-$r$ projector}}} 
   { \trace\big( \Sigma(I_d-P_r^T) H (I_d-P_r) \big)  }
  = { \sum_{i=r+1}^d \lambda_i } \, .
 \end{equation}
 Furthermore a solution to the above minimization problem is the $\Sigma^{-1}$-orthogonal projector defined by
 \begin{equation}\label{eq:OptimalPr}
  P_r = \Big(\sum_{i=1}^r v_iv_i^T \Big) \Sigma^{-1} \,.
 \end{equation}
\end{proposition}

By Propositions \ref{prop:BoundVVfunction} and \ref{prop:OptimalPr} we have that, for a sufficiently regular function $f$, the error $\| f - \mathbb{E}_\mu(f|\sigma(P_r)) \|_{\mathcal{H}}$ can be controlled by means of the generalized eigenvalues $\lambda_1,\hdots,\lambda_d$ of the matrix pair $(H,\Sigma^{-1})$ as follows
$$
 \| f - \mathbb{E}_\mu(f|\sigma(P_r)) \|_{\mathcal{H}}^2 \leq  { \sum_{i=r+1}^d \lambda_i }  ,
$$
where $P_r$ is the projector defined as in \eqref{eq:OptimalPr} and $H$ as in \eqref{eq:defH}. The matrix pair $(H,\Sigma^{-1})$ provides a test to reveal the low intrinsic dimension of the function $f$. Indeed, a fast decay in the spectrum of $(H,\Sigma^{-1})$ ensures that $\sum_{i=r+1}^d \lambda_i$ goes quickly to zero with $r$. In that case, given $\varepsilon>0$, there exists $r(\varepsilon) \ll d$ and a projector $P_r$ with rank $r(\varepsilon)$ such that $\| f - \mathbb{E}_\mu(f|\sigma(P_r)) \|_{\mathcal{H}}\leq\varepsilon$. Notice, however, that a fast decay in the spectrum of $(H,\Sigma^{-1})$ is only a \emph{sufficient} condition for the low intrinsic dimension: the absence of decay in the $(\lambda_i)$ does not mean that $f$ {cannot} be well approximated by $\mathbb{E}_\mu(f|\sigma(P_r))$ for some low-rank projector $P_r$.

\section{Contrast with the truncated Karhunen-Loève decomposition}\label{sec:KL}

A simple yet powerful dimension reduction method is the truncated Karhunen-Loève (K-L) decomposition. Truncated K-L decompositions are widely used in forward UQ---for instance, in any parameterized elliptic PDE where one truncates the K-L representation of the stochastic process modeling the uncertain parameters, and computes the resulting PDE solution; see, e.g., \cite{charrier2013finite,frauenfelder2005finite,nobile2016adaptive}.
The goal of this section is to position the method we proposed in Section \ref{sec2} against truncation of the K-L decomposition of the input parameters to $f$. In particular, under a Lipschitz continuity assumption on $f$, we show that K-L truncation can \textit{also} be seen as a function approximation technique that minimizes an upper bound. We then show that this upper bound is looser than the bound we derived using Poincar\'e inequalities.

The truncated K-L decomposition consists in reducing the parameter space to the subspace spanned by the leading eigenvectors of the covariance matrix of $\mu=\mathcal{N}(m,\Sigma)$. This approach is based on the observation that 
\begin{align}
 \min_{\substack{ P_r\in\mathbb{R}^{d\times d} \\ \text{rank-$r$ projector}}}
 &\mathbb{E}\big( \| (X-m) - P_r(X-m) \|_2^2 \big) \nonumber\\
 &=\min_{\substack{ P_r\in\mathbb{R}^{d\times d} \\ \text{rank-$r$ projector}}}
 \trace( (I_d-P_r)\Sigma(I_d-P_r^T))=\sum_{i=r+1}^d \sigma_i^2\, , \label{eq:KLmotivation}
\end{align}
where $X\sim\mu$ and where $\sigma_i^2$ is the $i$-th eigenvalue of $\Sigma$. 
We recall that $\|\cdot\|_2$ denotes the canonical norm of $\mathbb{R}^d$.
If the left-hand side of \eqref{eq:KLmotivation} is small, then the random variable $X$ can be well approximated (in the $L^2$ sense) by $m+P_r(X-m) = P_rX +(I_d-P_r)m$, where $P_r$ is a solution\footnote{Consider the eigendecomposition of $\Sigma = \sum_{i=1}^d \sigma^2_i u_iu_i^T$. Then the projector $P_r =\sum_{i=1}^r u_iu_i^T$ is a solution to \eqref{eq:KLmotivation}.} to \eqref{eq:KLmotivation}. In that case, given a function $f\in\mathcal{H}$, we can hope that $f(P_rX +(I_d-P_r)m)$ is a good approximation of $f(X)$. In order to make a quantitative statement, we assume $f$ is Lipschitz continuous, meaning that there exists a constant $L\geq0$ such that
\begin{equation}\label{eq:Lipschitz}
 \|f(x)-f(y)\|_V\leq L \| x-y \|_2 \,,
\end{equation}
for all $x,y\in\mathbb{R}^d$. Letting $g:x\mapsto f( P_r x +(I_d-P_r)m )$, we can write
\begin{align}
 \| f - g\circ P_r \|_{\mathcal{H}}^2
 &= \mathbb{E}\big( \| f(X) - f( P_rX +(I_d-P_r)m ) \|_V^2 \big) \nonumber\\
 &\overset{\eqref{eq:Lipschitz}}{\leq} L^2 \, \mathbb{E}\big( \| X - ( P_r X +(I_d-P_r)m ) \|_2^2 \big) 
 \overset{\eqref{eq:KLmotivation}}{=} L^2 \sum_{i=r+1}^d \sigma_i^2 . \label{eq:KL}
\end{align}
If the eigenvalues of $\Sigma$ decay rapidly, then there exist a function $g$ and a projector $P_r$ such that $\| f - g\circ P_r \|_{\mathcal{H}}\leq\varepsilon$, where $\text{rank}(P_r)=r(\varepsilon)\ll d$. In other words, the low intrinsic dimension of a Lipschitz continuous function can be revealed by the spectrum of $\Sigma$. Approximations that exploit this type of low-dimensional structure have been used extensively in forward and inverse uncertainty quantification; see, e.g., \cite{Marzouk2009}.

Notice that the function $g:x\mapsto f( P_r x +(I_d-P_r)m )$ considered here does not satisfy $g\circ P_r = \mathbb{E}_\mu(f|\sigma(P_r))$ in general, and therefore is not the optimal choice of profile; see Section \ref{sec:OptimalG}. 

\begin{proposition}\label{prop:KL_2}
 Let $f\in\mathcal{H}=L^2( \mathbb{R}^d , \mathcal{B}(\mathbb{R}^d) ,\mu ; V )$ be a continuously differentiable function and let $P_r$ be a minimizer of $P_r\mapsto \trace( \Sigma(I_d-P_r^T) H (I_d-P_r) ) $, where $H = \int (\nabla f)^T R_V \nabla f\mathrm{d}\mu$ and where $\Sigma$ is the covariance matrix of $\mu=\mathcal{N}(m,\Sigma)$. If $f$ is Lipschitz continuous such that \eqref{eq:Lipschitz} holds for some $L\geq0$, we have
 $$
  \| f - \mathbb{E}_\mu(f|\sigma(P_r)) \|_{\mathcal{H}}^2
  \leq \sum_{i=r+1}^d \lambda_i 
  \leq L^2  \sum_{i=r+1}^d \sigma_i^2 ,
 $$
 where $\sigma_i^2$ and $\lambda_i$ are the $i$-th eigenvalues of $\Sigma$ and of the matrix pair $(H,\Sigma^{-1})$ respectively.
 
\end{proposition}

The proof is given in Appendix \ref{proof:prop:KL_2}.
Similar to the methodology proposed in this paper, the truncated K-L decomposition can be interpreted as a method that minimizes an upper bound of an approximation error; see equation \eqref{eq:KL}. 
Proposition \ref{prop:KL_2} shows that the minimum of the upper bound of the new method is always smaller or equal to that of the truncated K-L. Of course comparing upper bounds does not allow one to make any clear statement about which method performs better than the other. However, note that for the truncated K-L decomposition, the construction of the projector relies only on the covariance matrix $\Sigma$, whereas the proposed method also takes into account the function $f$ (through the matrix $H$) in the construction of $P_r$. Thus it is natural to expect the new approach to provide projectors that are better for the approximation of $f$.

\section{Connection with global sensitivity measures}\label{sec:DGSM}

The goal of global sensitivity analysis is to assign,  to each group of input variables, a value that reflects its contribution to the variance of the output.
When considering a scalar-valued function $f:\mathbb{R}^d\rightarrow V$ with $V=\mathbb{R}$, classical variance-based indices include the closed Sobol' indices and the total Sobol' indices, defined respectively as: 
\begin{equation}\label{eq:defSobolS}
 S_\tau = \frac{\Var( \mathbb{E}(f(X)|X_\tau) ) }{\Var(f(X))} 
\quad\text{and}\quad
 T_\tau = 1-\frac{\Var( \mathbb{E}(f(X)|X_{-\tau}) ) }{\Var(f(X))} .
\end{equation}
Here $X_\tau$ and $X_{-\tau}$ represent components of the random vector $X\sim\mu$ indexed by $\tau$ and $-\tau$, where $\tau\subset\{1,\hdots,d\}$ is a set of indices with $\#\tau=r$ and where $-\tau$ is its complement in $\{1,\hdots,d\}$.
For independent inputs (e.g., diagonal $\Sigma$), the closed index $S_{\tau}$ measures $X_{\tau}$'s contribution to the output variance. The total index $T_\tau$ measures the contribution of $X_{\tau}$ and its interactions, of any order and with any other input variables, to the output variance.

The definitions in \eqref{eq:defSobolS} do not apply to vector-valued functions.
A natural extension of these indices is to interpret the variance of a (scalar-valued) function $h:\mathbb{R}^d\rightarrow \mathbb{R}$ as an $L^2$ norm, e.g., $\Var(h(X)) = \mathbb{E}(\|h(X)-\mathbb{E}(h(X))\|_V^2) $ where $V=\mathbb{R}$ with $\|\cdot\|_V=|\cdot|$. 
With this perspective, a natural extension of Sobol' indices to the vector-valued case $V\neq\mathbb{R}$ is
\begin{equation}\label{eq:defSobolS2}
 S_\tau 
 = \frac{\mathbb{E}( \| \mathbb{E}(f(X)|X_\tau) - \mathbb{E}(f(X)) \|_V^2)}{\mathbb{E}( \| f(X) - \mathbb{E}(f(X)) \|_V^2)}
 ~\text{and}~
 T_\tau 
 =  1-\frac{\mathbb{E}( \| \mathbb{E}(f(X)|X_{-\tau}) - \mathbb{E}(f(X)) \|_V^2)}{\mathbb{E}( \| f(X) - \mathbb{E}(f(X)) \|_V^2)} .
\end{equation}
Note that the definitions in \eqref{eq:defSobolS} and \eqref{eq:defSobolS2} are equivalent for scalar-valued functions.
We mention that a similar\footnote{To be specific, the generalization proposed in \cite{gamboa:hal-00800847,gamboa:hal-00881112} is $S_\tau=\trace(M C_\tau)/\trace(MC)$, where $C_\tau=\text{Cov}(\mathbb{E}(f(X)|X_\tau))$ and $C=\text{Cov}(f(X))$, where $M\in\mathbb{R}^{n\times n}$ is a given matrix. One can easily show that if $M=R_V$, which means $\|y\|_V^2=y^TMy$ for any $y\in\mathbb{R}^n$, then this definition matches the one proposed in \eqref{eq:defSobolS2}.} generalization of the Sobol' index $S_\tau$ has been proposed in \cite{gamboa:hal-00800847,gamboa:hal-00881112}.
Using standard properties of the conditional expectation, one can rewrite the above indices as 
$$
 S_\tau 
 = 1-\frac{\mathbb{E}( \| f(X)- \mathbb{E}(f(X)|X_\tau) \|_V^2)}{\mathbb{E}( \| f(X) - \mathbb{E}(f(X)) \|_V^2)}
 = 1-\frac{\| f-\mathbb{E}_\mu(f|\sigma(P_r)) \|_\mathcal{H}^2}{\| f- \mathbb{E}_\mu(f) \|_\mathcal{H}^2} ,
$$
and
$$
 T_\tau 
 =  \frac{\mathbb{E}( \| f(X)-\mathbb{E}(f(X)|X_{-\tau})  \|_V^2)}{\mathbb{E}( \| f(X) - \mathbb{E}(f(X)) \|_V^2)} 
 = \frac{\| f-\mathbb{E}_\mu(f|\sigma(I_d-P_r)) \|_\mathcal{H}^2}{\| f- \mathbb{E}_\mu(f) \|_\mathcal{H}^2},
$$
where $P_r$ is the projector such that $P_rX$ (resp.\ $(I_d-P_r)X$) extracts the coordinates of $X$ indexed by $\tau$ (resp.\ by $-\tau$).
As noticed in \cite{hartapproximation}, the above expressions allow for an interpretation of the Sobol' indices with an approximation perspective. On the one hand, $S_\tau$ quantifies how well a function $f$ can be approximated by $\mathbb{E}_\mu(f|\sigma(P_r))$, a function which depends only on the $\tau$-coordinates of the parameter (large $S_\tau$ means we should not remove the $X_\tau$ dependence).
On the other hand, $T_\tau$ quantifies how good an approximation of $f$ can be if we remove the coordinates indexed by $\tau$ (small $T_\tau$ means we can remove the $X_\tau$ dependence). 

A straightforward application of Proposition \ref{prop:BoundVVfunction} allows us to bound the indices $S_\tau$ and $T_\tau$ as follows:
\begin{equation}\label{eq:lowerboundS}
 S_\tau 
 \geq 1-\frac{\trace ( \Sigma(I_d-P_r^T) H (I_d-P_r) )}{\| f- \mathbb{E}_\mu(f) \|_\mathcal{H}^2}
 = 1-\frac{\sum_{i\notin\tau} \Var(X_i)H_{i,i} }{\| f- \mathbb{E}_\mu(f) \|_\mathcal{H}^2},
\end{equation}
and
\begin{equation}\label{eq:upperboundT}
 T_\tau 
 \leq \frac{\trace ( \Sigma(P_r^T) H (P_r) )}{\| f- \mathbb{E}_\mu(f) \|_\mathcal{H}^2}
 = \frac{\sum_{i\in\tau} \Var(X_i) H_{i,i} }{\| f- \mathbb{E}_\mu(f) \|_\mathcal{H}^2} ,
\end{equation}
where $\Var(X_i) = \Sigma_{i,i}$.
In the scalar-valued case, $H_{i,i} = \int (\partial_{i}f)^2\mathrm{d}\mu$ coincides with the $i$th derivative-based global sensitivity measure (DGSM) \cite{kucherenko2009derivative,kucherenko2009monte}. 
The fact that the DGSM can bound the total Sobol' index $T_\tau$ has been already noted \cite{kucherenko2009derivative,kucherenko2016derivative,lamboni2013derivative,roustant2017poincare} for scalar-valued functions, and for more general input distributions than Gaussian ones.
Here, under the assumption of Gaussian probability measure $\mu$, inequality \eqref{eq:upperboundT} provides a generalization of these bounds to the case of vector-valued functions, where the $i$th DGSM ought to be defined as $H_{i,i} = \int \|\partial_{i}f\|_V^2\mathrm{d}\mu$.
The same remark applies for inequality \eqref{eq:lowerboundS}.

\section{Illustrations}\label{sec4}

\subsection{Analytical examples}
\label{sec:analytical}

We give here three analytical examples for which we can compute a closed-form expression for the error $\| f - \mathbb{E}_\mu(f|\sigma(P_r)) \|_{\mathcal{H}}$. This allows us to find the projector that minimizes the true error. We then compare this projector with the one that minimizes the upper bound of $\| f - \mathbb{E}_\mu(f|\sigma(P_r)) \|_{\mathcal{H}}$.

First we consider a linear function. We show that the bound equals the true error, so that minimizing the bound gives the minimizer of the error itself. Then we consider a quadratic function: in this case, the bound is not equal to the error, but the minimizers are the same.
Finally we consider a function defined as a sum of sine functions. Depending on the frequency and amplitude of the sines, minimizing the bound can either yield the optimal projector (minimizer of the error) or the worst projector (maximizer of the error)! This last example offers some useful intuition, showing that the proposed method performs better for slowly varying functions than for functions of small amplitude but high frequency.

\subsubsection{Linear functions}

Assume $f\in\mathcal{H}$ is a linear function $f:x\mapsto Fx$ for some matrix $F\in\mathbb{R}^{n\times d}$ and let $P_r\in\mathbb{R}^{d\times d}$ be a $\Sigma^{-1}$-orthogonal projector. 
By Proposition \ref{prop:ExplicitCondExp} and by linearity of $f$ we have
$\mathbb{E}_\mu(f|\sigma(P_r))(x) = FP_rx+ F(I_d-P_r)m$ for any $x\in\mathbb{R}^d$. We can write
\begin{align*}
\| f- \mathbb{E}_\mu(f|\sigma(P_r))\|_\mathcal{H}^2 
&= \int_{\mathbb{R}^d} \| Fx - FP_rx - F(I_d-P_r)m \|_V^2 \,\mathrm{d}\mu(x) \\
&= \int_{\mathbb{R}^d} \| F(I_d-P_r)(x - m) \|_V^2 \,\mathrm{d}\mu(x) \\
&= \int_{\mathbb{R}^d} (x - m)^T(I_d-P_r)^T F^T R_V F (I_d-P_r)(x - m)  \,\mathrm{d}\mu(x) \\
&= \text{trace} \big( \Sigma (I_d-P_r^T) H (I_d-P_r) \big)  \, ,
\end{align*}
where, for the last equality, we used the relations $\Sigma = \int_{\mathbb{R}^d} (x - m)(x - m)^T\mathrm{d}\mu(x)$ and $H = \int (\nabla f)^TR_V(\nabla f) \mathrm{d}\mu =  F^TR_VF$.
Thus we have that equality is attained in \eqref{eq:BoundVVfunction} for any linear functions $f\in\mathcal{H}$ and for any $\Sigma^{-1}$-orthogonal projector $P_r$. This shows that, for linear functions, the upper bound is equal to the true error.

\subsubsection{Quadratic forms}

Assume $\mu=\mathcal{N}(0,I_d)$ is the standard normal distribution and let $f\in\mathcal{H}$ be a quadratic form defined by $f:x\mapsto \frac{1}{2}x^TAx$ for some symmetric matrix $A\in\mathbb{R}^{d\times d}$. It is a real-valued function so that $V=\mathbb{R}$ and $\|\cdot\|_V = |\cdot|$, the absolute value. 
Let $P_r$ be an orthogonal projector with rank $r$ so that $P_r^T=P_r$. One can easily check that the relation 
$$
f(P_rx+(I_d-P_r)Y) = f(P_rx) + Y^T(I_d-P_r)AP_rx + f((I_d-P_r)Y) \,,
$$
holds for all $x\in\mathbb{R}^d$ where $Y\sim\mu$. By taking the expectation with respect to $Y$, Proposition \ref{prop:ExplicitCondExp} allows writing
$
\mathbb{E}_\mu(f|\sigma(P_r))(x) = f(P_rx) + \mathbb{E}( f((I_d-P_r)Y) ).
$
The function $f-\mathbb{E}_\mu(f|\sigma(P_r))$ is quadratic and can be written as $x\mapsto x^T \Lambda x + c$ where $\Lambda = \frac{1}{2}(A - P_rAP_r)$ and $c = - \mathbb{E}( Y^T\Lambda Y )$. 
We have
\begin{align*}
 \| f - \mathbb{E}_\mu(f|\sigma(P_r)) \|_{\mathcal{H}}^2 
 = \mathbb{E}\big( (Y^T\Lambda Y + c)^2 \big)
 = \Var( Y^T \Lambda Y ) \,. 
\end{align*}
Consider the eigendecomposition of $\Lambda=U\text{diag}(a_1,\hdots,a_d)U^T$ and let $Z=U^TY \sim \mathcal{N}(0,I_d)$. We have $Y^T \Lambda Y= \sum_{i=1}^d a_i Z_i^2 $ so that
\begin{align*}
 \| f - \mathbb{E}_\mu(f|\sigma(P_r)) \|_{\mathcal{H}}^2 
 = \sum_{i=1}^d a_i^2 \Var(  Z_i^2 ) = 2\sum_{i=1}^d a_i^2
 = 2\trace( \Lambda^2 )  = \frac{1}{2} \| A - P_rAP_r \|_F^2 \,,
\end{align*}
where $\|\cdot\|_F=\sqrt{\trace(\cdot)^T(\cdot)}$ denotes the Frobenius norm. One can show that the rank-$r$ projector which minimizes $P_r\mapsto \| A-P_rAP_r \|_F$ is the projector onto the leading eigenspace of $A^2$. Denoting by $\alpha_i^2$ the $i$-th largest eigenvalue of $A^2$, we have
\begin{equation}\label{eq:tmp65261}
\min_{\substack{ P_r\in\mathbb{R}^{d\times d} \\ \text{rank-$r$ orth. projector}}} \| f - \mathbb{E}_\mu(f|\sigma(P_r)) \|_{\mathcal{H}} 
= \frac{1}{\sqrt{2}} \Big( \sum_{i=r+1}^d \alpha_i^2 \Big)^{1/2}.
\end{equation}

Now we consider the projector that minimizes the upper bound given by Proposition \ref{prop:BoundVVfunction}. We can write $\nabla f(x) = Ax$ so that $H = \int (\nabla f)(\nabla f)^T\mathrm{d}\mu =  A^2 $. Therefore equation \eqref{eq:BoundVVfunction} yields
$$
\| f - \mathbb{E}_\mu(f|\sigma(P_r)) \|_{\mathcal{H}}^2 
\leq \trace\big( (I_r-P_r) A^2 (I_r-P_r) \big)
= \|A-P_rA\|_F^2
\,,
$$
for any orthogonal projector $P_r$ with rank $r$. By Proposition \ref{prop:OptimalPr}, the rank-$r$ orthogonal projector which minimizes the right-hand side in the above inequality is the projector onto the leading eigenspace of $A^2$, which is the same as the solution to \eqref{eq:tmp65261}. Then the minimizer of the bound is, for the considered example, the same as the minimizer of the error itself. In addition, the upper bound evaluated at the optimal projector allows controlling the error $\| f - \mathbb{E}_\mu(f|\sigma(P_r)) \|_{\mathcal{H}}$ by $(\sum_{i>r}^d \alpha_i^2)^{1/2}$ which is, up to a factor of $\sqrt{2}$, the same as the true error.

\subsubsection{Sum of sines}

Let $\mu=\mathcal{N}(0,I_d)$ be a standard normal distribution. Consider the real-valued function $f\in\mathcal{H}$ such that
$$
f : x\mapsto \sum_{i=1}^d a_i \sin( \omega_i x_i ) \,,
$$
for any $x\in\mathbb{R}^d$, where $a\in\mathbb{R}^d$ and $\omega\in\mathbb{R}^d$ are two vectors. Let $P_r$ be an orthogonal projector. For simplicity, we restrict our analysis to the case where $P_r$ is a projector onto the span of $r$ vectors from the canonical basis $\{e_1,\hdots, e_d\}$ of $\mathbb{R}^d$, meaning 
\begin{equation}\label{eq:tmp23608}
P_r = \sum_{i\in \tau} e_i e_i^T, 
\end{equation}
where $\tau\subset\{1,\hdots,d\}$ and $\# \tau = r$.
It is readily seen that $\mathbb{E}_\mu(f|\sigma(P_r))$ is the function $ x\mapsto \sum_{i\in\tau} a_i \sin( \omega_i x_i )$. We can show that
$$
\| f - \mathbb{E}_\mu(f|\sigma(P_r)) \|_{\mathcal{H}}^2 
= \mathbb{E} \Big( \big( \sum_{i\in-\tau} a_i \sin( \omega_i X_i ) \big)^2 \Big)
= \frac{1}{2}\sum_{i\in-\tau} a_i^2 (1-\exp(-2\omega_i^2)) \,,
$$
where $-\tau$ is the complementary set of $\tau$ in $\{1,\hdots,d\}$ and $X\sim\mu$. Therefore, the projector $P_r$ of the form of \eqref{eq:tmp23608} which minimizes the error $\| f - \mathbb{E}_\mu(f|\sigma(P_r)) \|_{\mathcal{H}}$ is the one associated with the set $\tau$ containing the indices of the $r$ largest values of $a_i^2 (1-\exp(-2\omega_i^2))$.
\\

Now we find the projector of the form \eqref{eq:tmp23608} that minimizes the upper bound of the error given by Proposition \ref{prop:BoundVVfunction}. Recall that $H = \int (\nabla f)(\nabla f)^T\mathrm{d}\mu$, so we can write 
\begin{align*}
\| f - \mathbb{E}_\mu(f|\sigma(P_r)) \|_{\mathcal{H}}^2 
&\overset{\eqref{eq:BoundVVfunction}}{\leq} \trace( (I_d-P_r)^TH(I_d-P_r) ) \\
&\overset{\eqref{eq:tmp23608}}{=} \sum_{i\in-\tau} e_i^T H e_i = \sum_{i\in-\tau}  \int \Big(\frac{\partial f}{\partial x_i}\Big)^2 \mathrm{d}\mu \\
&= \sum_{i\in-\tau}   \mathbb{E}\big( ( a_i\omega_i \cos(\omega_i X_i) )^2 \big) \\
&= \frac{1}{2}\sum_{i\in-\tau}  a_i^2\omega_i^2 ( 1 + \exp(-2\omega_i^2) ) \,.
\end{align*}
The projector \eqref{eq:tmp23608} that minimizes the above upper bound is the one associated with the set $\tau$ containing the indices of the $r$ largest values of $a_i^2\omega_i^2 ( 1 + \exp(-2\omega_i^2) )$. 
We now describe two interesting cases.
\begin{itemize}
 \item Assume that all the frequencies are the same, $\omega_i = \omega$ for all $i\leq d$. The index sets corresponding to the largest $a_i^2\omega^2 ( 1 + \exp(-2\omega^2) )$ and $a_i^2 (1-\exp(-2\omega^2))$ are the same, and therefore the projector that minimizes the upper bound is the same as the minimizer of the true error. Notice, however, that when $\omega\rightarrow \infty$ the true error tends to $\frac{1}{2}\sum_{i\in-\tau} a_i^2$, whereas the upper bound tends to infinity. This shows that the upper bound can be a poor estimator for the error, even if its minimization allows recovery of the optimal projector.

 \item Suppose now that $\omega_i = a_i^{-2} \geq 1$ for all $i\leq d$. Then the index set corresponding to the largest $a_i^2\omega_i^2 ( 1 + \exp(-2\omega_i^2) ) =  \omega_i(1 + \exp(-2\omega_i^2)) \eqqcolon h_1(\omega_i)$ is the same as the index set of the smallest $a_i^2 (1-\exp(-2\omega_i^{2}))=\omega_i^{-1} (1-\exp(-2\omega_i^{2})) \eqqcolon h_2(\omega_i)$. Indeed $h_1$ is increasing on $(1,\infty)$ whereas $h_2$ is decreasing. Hence, for this particular example, minimizing the upper bound yields the worst possible projector, i.e., the one that maximizes the true error.
 
\end{itemize}

These two cases show the limitations of the use of Poincaré inequalities: the bound is not sharp for functions with small variation but high frequencies. However, it works well for slowly varying functions.
The same remark applies directly to sensitivity analysis (see Section \ref{sec:DGSM}): the DGSM should not be used to bound the Sobol' indices unless the function varies slowly with respect to its input parameters.

\subsection{Elliptic PDE}
\label{sec:pde}

\begin{figure}
    \centering
    \begin{subfigure}[b]{0.31\textwidth}\centering
        \includegraphics[width=0.95\textwidth]{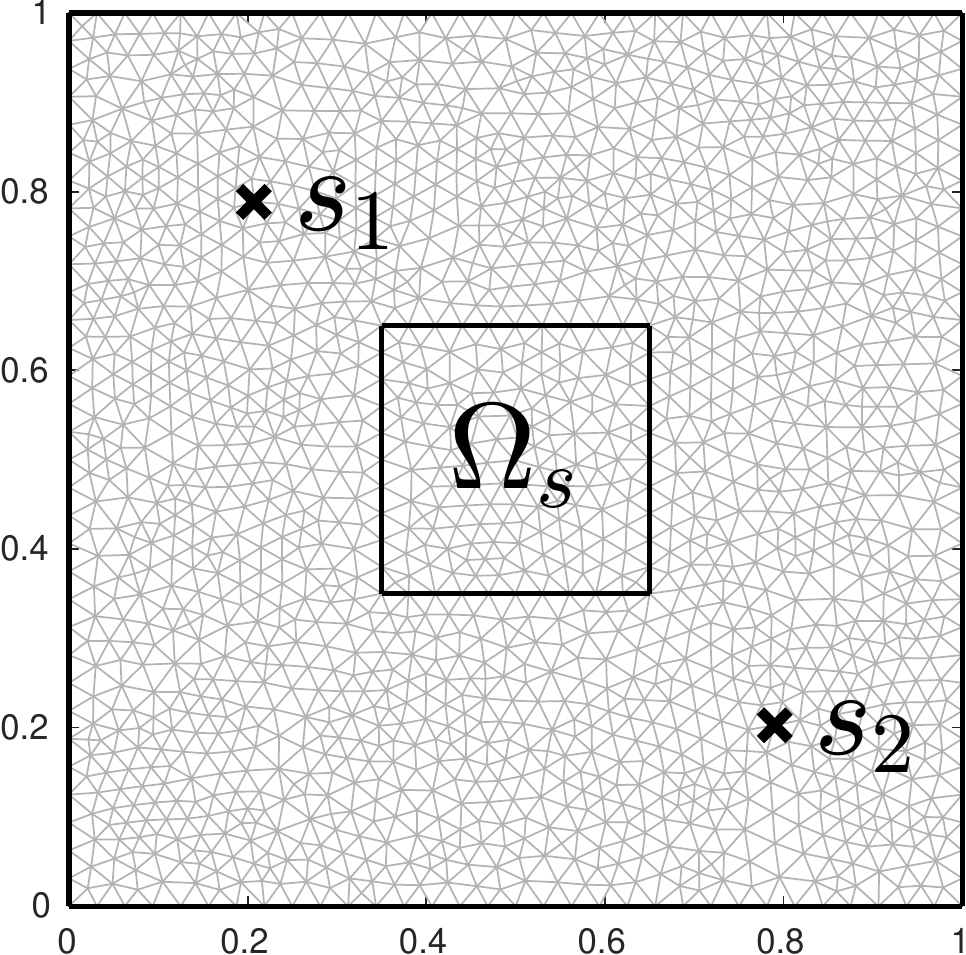}
        \caption{Geometry and mesh.}
        \label{fig:domain}
    \end{subfigure}
    ~
    \begin{subfigure}[b]{0.31\textwidth}\centering
        \includegraphics[width=0.96\textwidth]{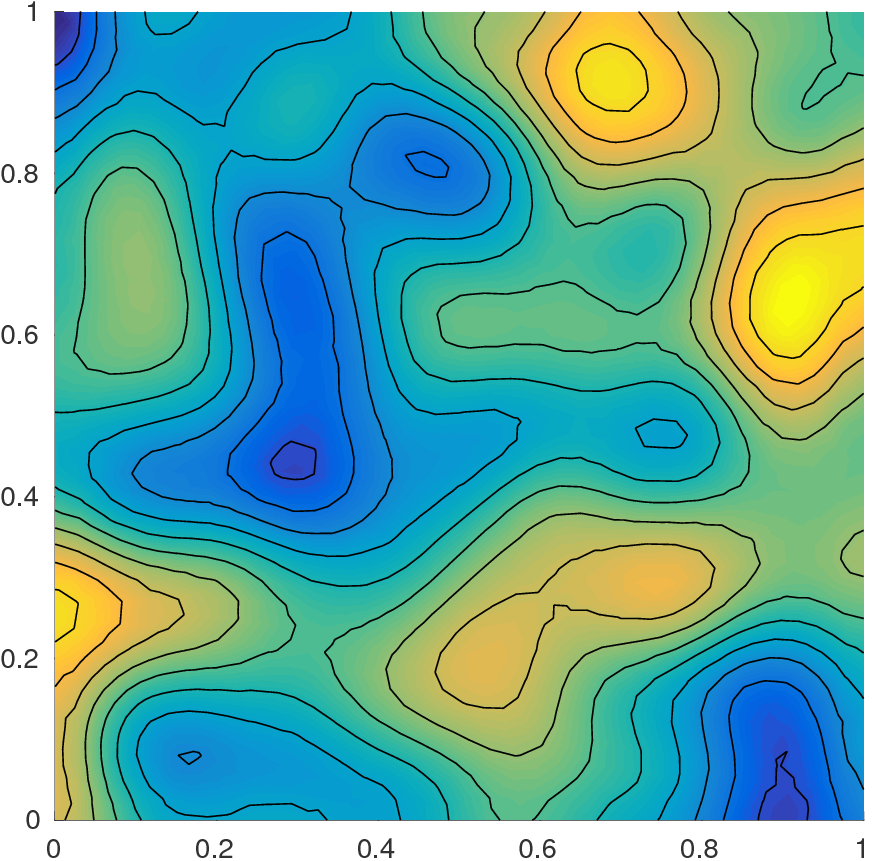}
        \caption{Diffusion field $\log(\kappa(x))$.}
        \label{fig:diffusionField}
    \end{subfigure}
    ~
    \begin{subfigure}[b]{0.31\textwidth}\centering
        \includegraphics[width=0.96\textwidth]{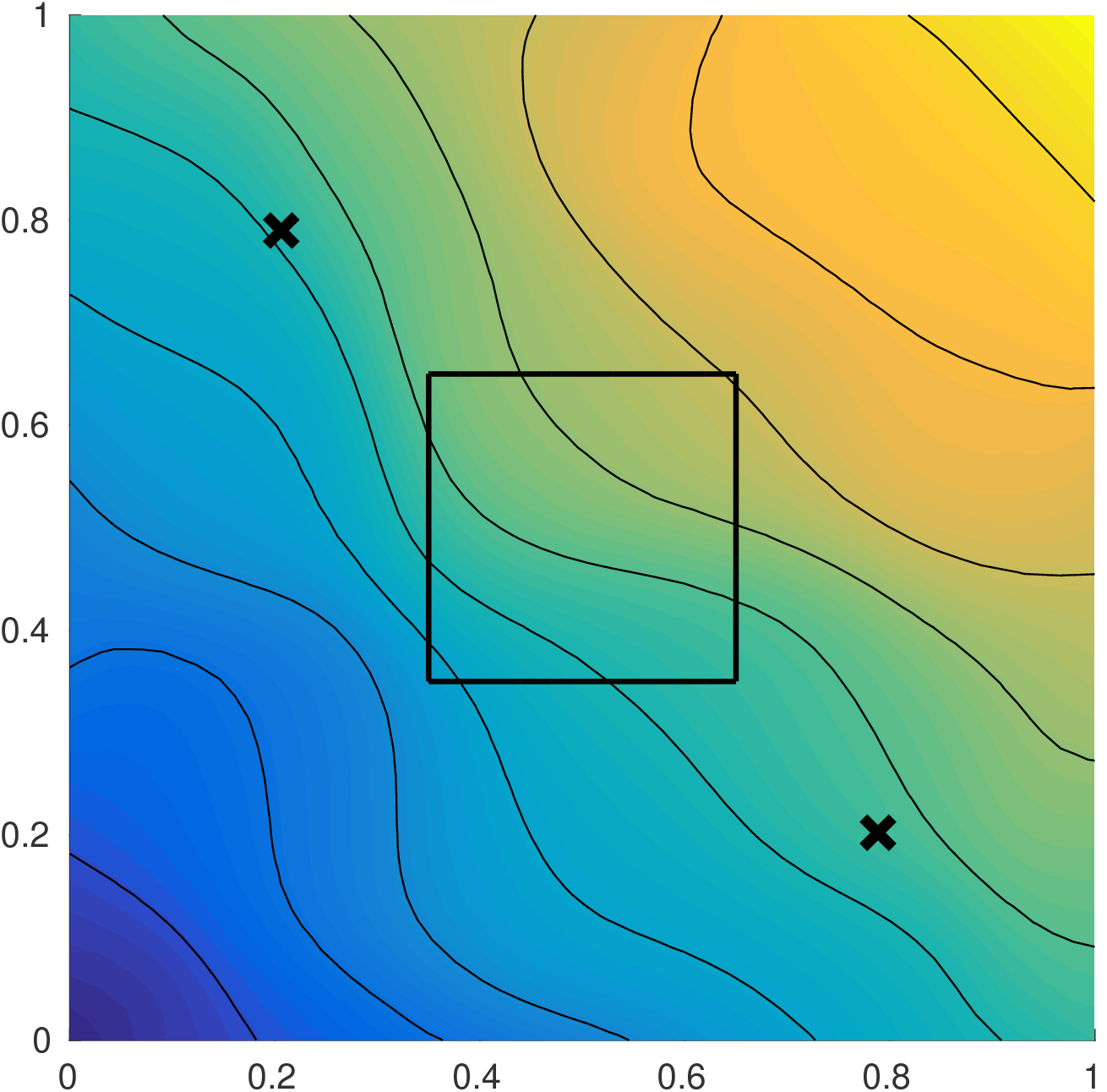}
        \caption{Solution $u(x)$.}
        \label{fig:solution}
    \end{subfigure} 
    \\[0.2cm]
    \begin{subfigure}[b]{0.9\textwidth}\centering
      \begin{tabular}{|cm{0.2\textwidth}|cm{0.2\textwidth}|m{0.21\textwidth}|} \hline 
       \multicolumn{2}{|c|}{Scenario 1} & \multicolumn{2}{c|}{Scenario 2} & \multicolumn{1}{c|}{Scenario 3} \\ \hline
       $f(x)=$ & \vspace{0.1cm}\includegraphics[width=0.19\textwidth]{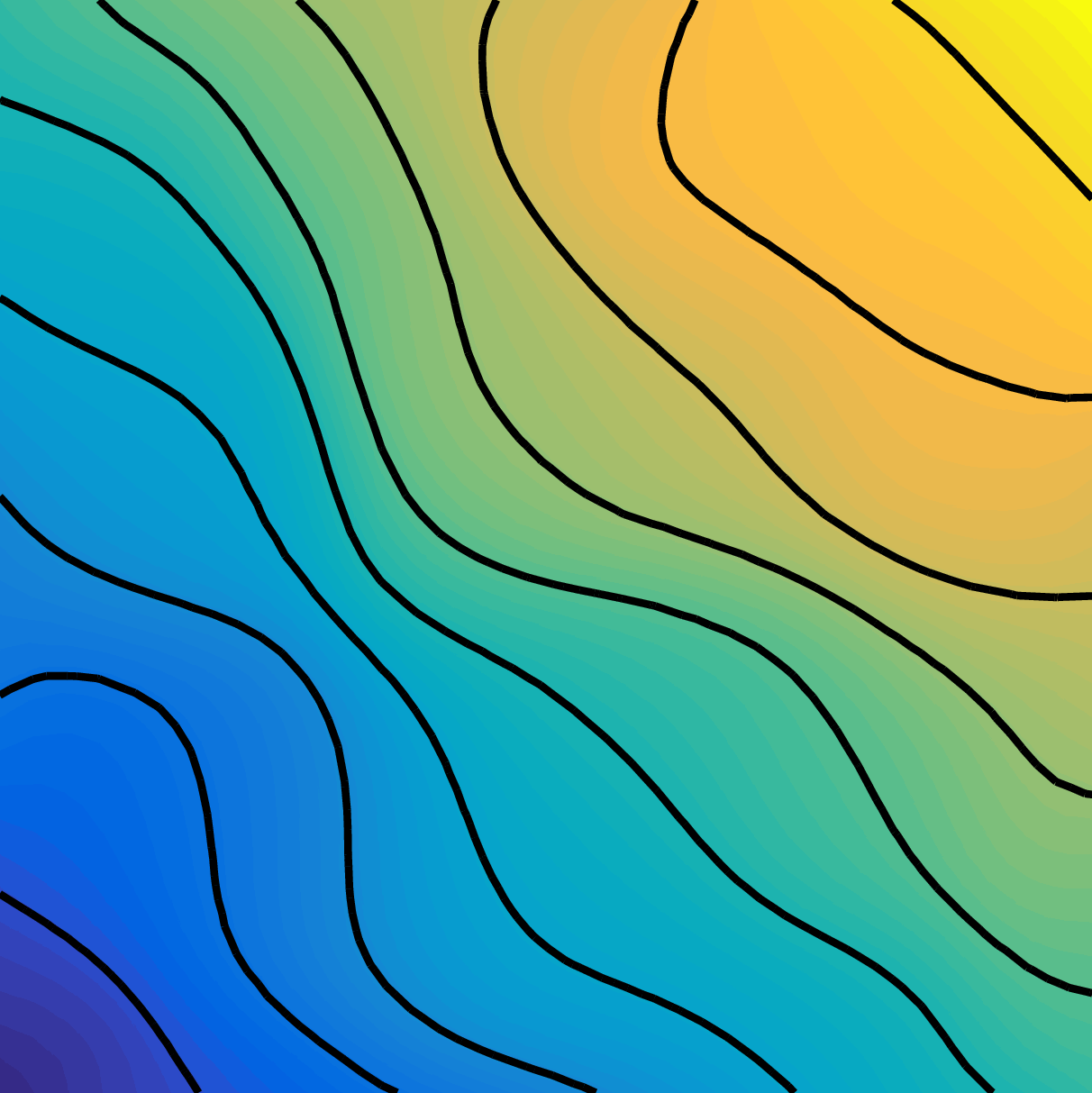} &
       ~~$f(x)=$ & \includegraphics[width=0.09\textwidth]{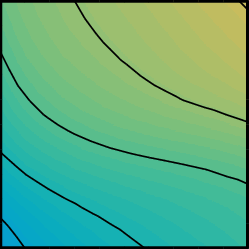} &
       $$ f(x)=\begin{pmatrix} 0.939 \\ 1.032 \end{pmatrix} $$ \\\hline
      \end{tabular}
      \caption{Three different scenarios.}
      \label{fig:postT}
    \end{subfigure}
    
  \caption{Illustration of the elliptic PDE problem: geometry and mesh (Figure \ref{fig:domain}), representation of the diffusion field associated with a parameter $x\in\mathbb{R}^d$ drawn randomly from $\mu$ (Figure \ref{fig:diffusionField}), corresponding solution (Figure \ref{fig:solution}) and representation of $f(x)$ for the three different scenarios given this particular $x$ (Figure \ref{fig:postT}).}
  \label{fig:EllipticPDE}
\end{figure}

Consider the diffusion equation on the square domain $\Omega=[0,1]^2$, which consists in finding $u$ in the Sobolev space $H^1(\Omega)$ such that
\begin{equation}\label{eq:diffProblem}
 \left\{
 \begin{array}{rl}
  \nabla_s (\kappa \, \nabla_s u ) &=0 ~~\quad\quad\quad \text{in } \Omega \,, \\
  u &=s_1+s_2 \quad \text{on } \partial\Omega \,.\\
 \end{array}
 \right.
\end{equation}
Here $s=(s_1,s_2)\in\Omega$ denotes the spatial coordinates and $\nabla_s$ refers to the gradient in the spatial variable $s$. 
The diffusion coefficient $\kappa$ is a random field and follows a log-normal distribution such that $\log(\kappa)$ is a Gaussian process on $\Omega$ with zero mean and with a covariance function $c:\Omega\times\Omega\rightarrow\mathbb{R}$ defined by
$
 c(s,t) = \exp ( - \|s-t\|_2^2 / (0.15)^2 ) 
$
for all $s,t\in\Omega$. 
A numerical approximation of \eqref{eq:diffProblem} is obtained with the finite element method (FEM); see, for example, \cite{Ern2004}. The diffusion field $\kappa$ is approximated by the piecewise constant random field
\begin{equation}\label{eq:defKappa}
 \kappa(x) :s\mapsto \exp \Big( \sum_{i=1}^{d} x_i \, \mathbf{1}_i(s)  \Big)\, ,
\end{equation}
where $\mathbf{1}_i$ denotes the indicator function associated with the $i$th element of the mesh represented in Figure \ref{fig:domain}. Here $d=3252$ corresponds to the number of elements, and $x\sim\mu=\mathcal{N}(0,\Sigma)$ with
$$
 \Sigma_{i,j} = c(s_i,s_j)
 ,\quad 1\leq i,j \leq d \,,
$$
and $s_i$ being the center of the $i$th element. With a slight abuse of notation, we denote by $u(x)$ the Galerkin projection of the solution to \eqref{eq:diffProblem} onto the space of continuous piecewise affine functions associated with the mesh in Figure \ref{fig:domain}. We consider the following scenarios, where the function $f:\mathbb{R}^d\rightarrow V$ is defined by three different post-solution treatments of $u(x)$:

\begin{enumerate}
 \item $f:x\mapsto u(x)$, which means that $f$ is the solution map from the parameter $x$ to the FEM solution to \eqref{eq:diffProblem}. In that case $V$ is the FEM approximation space with dimension $\text{dim}(V)=n=1691$, the number of nodes in the mesh. Since $V\subset H^1(\Omega)$, the natural choice for the norm $\|\cdot\|_V$ is
 $$
  \|v\|_V^2 = \int_\Omega (v(s))^2 \,\mathrm{d}s + \int_\Omega \|\nabla_s v(s)\|_2^2 \, \mathrm{d}s \,.
 $$
 
 \item $f:x\mapsto u_{|\Omega_s}(x)$, where $\Omega_s=[0.35 , 0.65]^2\subset\Omega$. In other words, $f(x)$ corresponds to the restriction of $u(x)$ to a subdomain $\Omega_s$ of $\Omega$. For this scenario, $V\subset H^1(\Omega_s)$ is of dimension $n=168$ (the number of nodes in $\Omega_s$) and is endowed with the norm $\|\cdot\|_V$ given by
 $$
  \|v\|_V^2 = \int_{\Omega_s} (v(s))^2 \,\mathrm{d}s + \int_{\Omega_s} \|\nabla_s v(s)\|_2^2 \, \mathrm{d}s \,.
 $$
 
 \item $f:x\mapsto (u_{|s_a}(x),u_{|s_b}(x))$, where $s_a=(0.2,0.8)\in\Omega$ and $s_b=(0.8,0.2)\in\Omega$. In this scenario, we are interested in the evaluation of the solution $u(x)$ at two different spatial locations $s_a$ and $s_b$. There are two scalar-valued outputs so that $V = \mathbb{R}^2$ is an algebraic space. Consider the weighted norm $\|\cdot\|_V$ defined by
 $$
  \|v\|_V^2 = \alpha \, v_1^2 +  \beta \, v_2^2 \,,
 $$
 where $\alpha,\beta>0$ are two positive weights to be specified.
 For example the choice $\alpha = 2\beta$ will put twice the weight on the error associated with the first output compared to the second. This is a way to model the fact that, for the final purpose of the simulation, one output is more important than the other.
 \\
 
\end{enumerate}

\begin{remark}
 Each of the three functions defined above are continuously differentiable, as a composition of continuously differentiable functions. Indeed, we can write each $f:x\mapsto Lu(x)$ with some matrix $L$ which depends on the scenario and $u(x)=A(x)^{-1}b$ where $A(x)$ and $b$ are the FEM matrix and FEM right-hand side associated with \eqref{eq:diffProblem}. 
 Because of the parameterization \eqref{eq:defKappa} of the diffusion field $\kappa(x)$, the function $x\mapsto A(x)$ is continuously differentiable so that $x\mapsto u(x)=A(x)^{-1}b$ is also continuously differentiable, and so is $f:x\mapsto Lu(x)$.
\end{remark}

\subsubsection{Computational aspects}
\label{sec:compuational}
We consider the problem of computing the matrix $H=\int (\nabla f)^TR_V ( \nabla f)\, \mathrm{d} \mu$. Since $H=\mathbb{E}( (\nabla f(X))^TR_V ( \nabla f(X)) )$, with $X\sim \mu$, $H$ can be approximated by the $K$-sample Monte Carlo estimate
\begin{equation}\label{eq:defHhat}
 \widehat H = \frac{1}{K} \sum_{i=1}^K \big(\nabla f(X_i) \big)^TR_V \big( \nabla f(X_i) \big) \,,
\end{equation}
where $X_1,\hdots,X_K$ are independent copies of $X$. To numerically compute a realization of $\widehat H$, one needs to evaluate the Jacobian of the function $f$ $K$ times. To do so, we employ the \emph{adjoint method}; see for example \cite{plessix2006review}.
Then, to construct the projector, instead of minimizing $\trace( \Sigma(I_d-P_r^T) H (I_d-P_r) )$ we consider a projector $\widehat P_r$ such that
\begin{equation}\label{eq:defPrHat}
 \widehat P_r\in\underset{\substack{ P_r\in\mathbb{R}^{d\times d} \\ \text{rank-$r$ projector}}} {\text{arg} \> \text{min}}
 \trace\big( \Sigma(I_d-P_r^T) \widehat H (I_d-P_r) \big) \,.
\end{equation}
By construction, $\widehat P_r$ depends upon $\widehat H$, and thus it is random. Recall that such a projector can be obtained by computing the generalized eigendecomposition of the matrix pair $(\widehat H,\Sigma^{-1})$; see Proposition \ref{prop:OptimalPr}.

To approximate the conditional expectation $\mathbb{E}_\mu(f|\sigma(\widehat P_r))$, we consider the random function
\begin{equation}\label{eq:defFhatr}
 \widehat F_r : x\mapsto \frac{1}{M} \sum_{i=1}^M  f(\widehat P_r x + (I_d-\widehat P_r)Y_i) \,,
\end{equation}
where $Y_1,\hdots,Y_M$ are independent copies of $Y\sim\mu$.  
Given a realization of the projector $\widehat P_r$, a realization of $\widehat F_r$ can be obtained by drawing $M$ samples of $Y$ and by using those samples to evaluate $\widehat F_r$ using \eqref{eq:defFhatr}. Notice that the samples are \emph{not} redrawn for each new evaluation point $x$ of $\widehat F_r$. By Proposition \ref{prop:ExplicitCondExp} and for any $x\in\mathbb{R}^d$, $\widehat F_r(x)$ can be interpreted as an $M$-sample Monte Carlo approximation of $\mathbb{E}_\mu(f|\sigma(\widehat P_r))(x)$. Finally, notice that if $M=1$ and $Y_1 = 0$ (i.e., the mean of $Y$), then our approximation of $f$ reduces to the form used in Section~\ref{sec:KL} when truncating a K-L decomposition, albeit for a different projector; see relation \eqref{eq:KL} with $m=0$ and $P_r=\widehat P_r$.

\subsubsection{Modes and influence of the norm $\|\cdot\|_V$}

\begin{figure}\centering
\begin{tabular}{ccccccc}
 & mode 1 & mode 2 & mode 3 & mode 4 & mode 5 & mode 6\\
\rotatebox[origin=l]{90}{~~~~K-L}&
\includegraphics[width=0.115\textwidth]{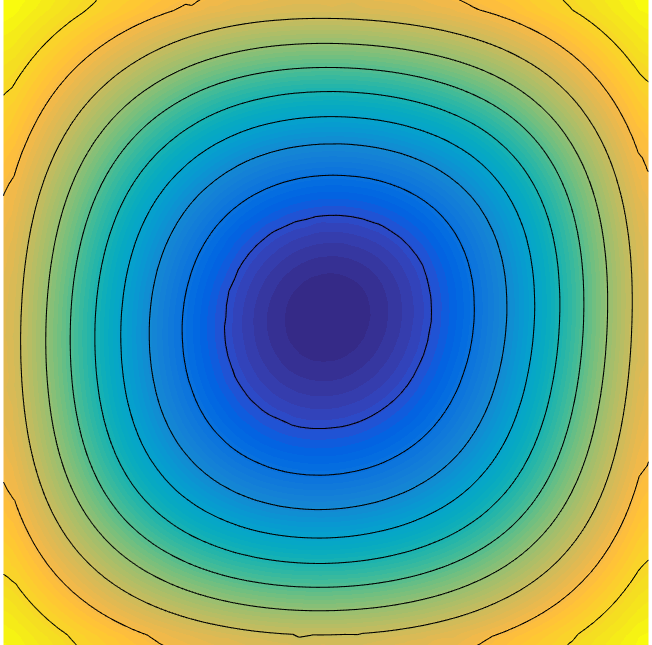} & 
\includegraphics[width=0.115\textwidth]{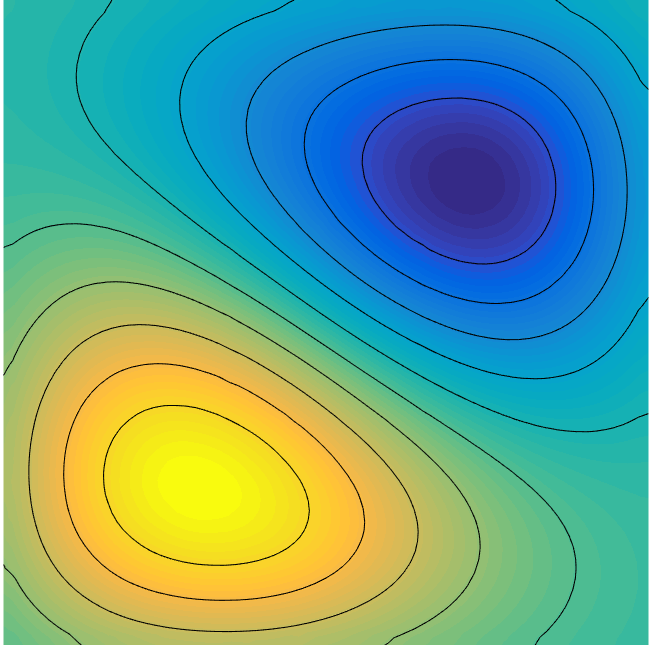} & 
\includegraphics[width=0.115\textwidth]{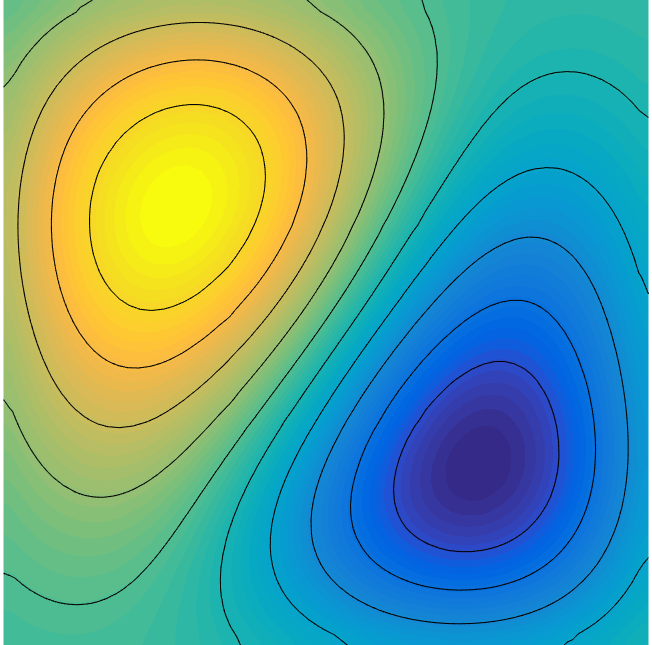} & 
\includegraphics[width=0.115\textwidth]{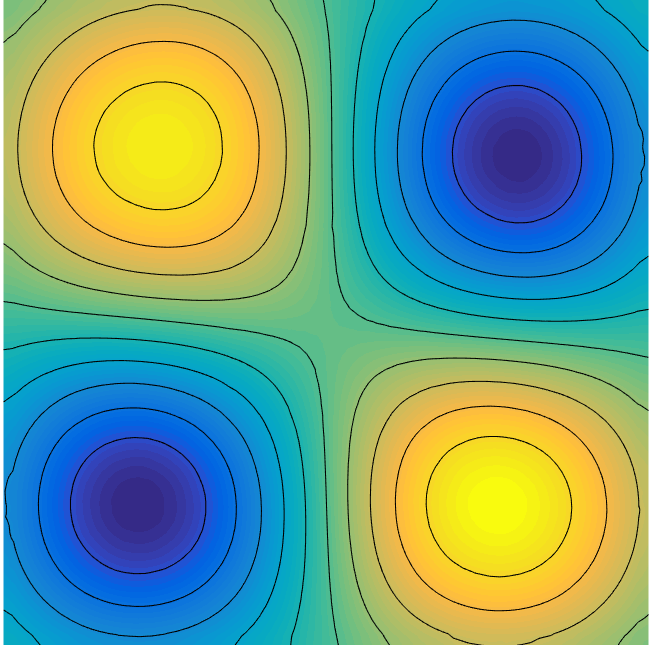} & 
\includegraphics[width=0.115\textwidth]{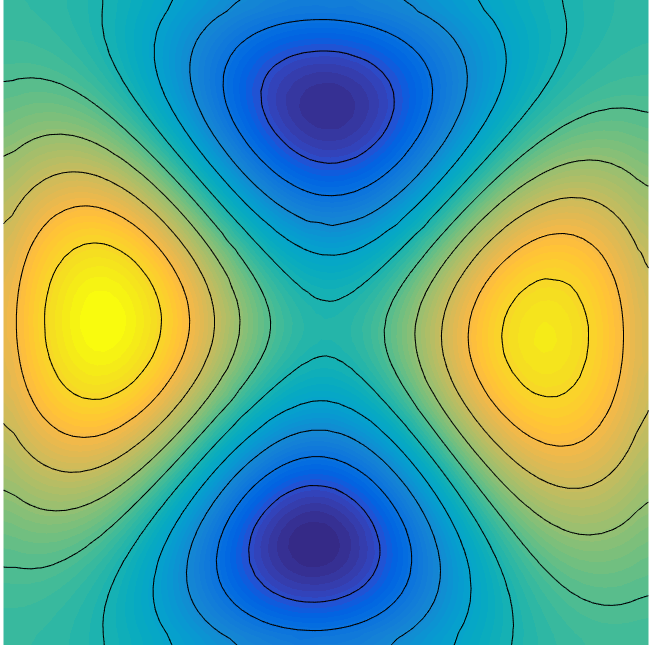} & 
\includegraphics[width=0.115\textwidth]{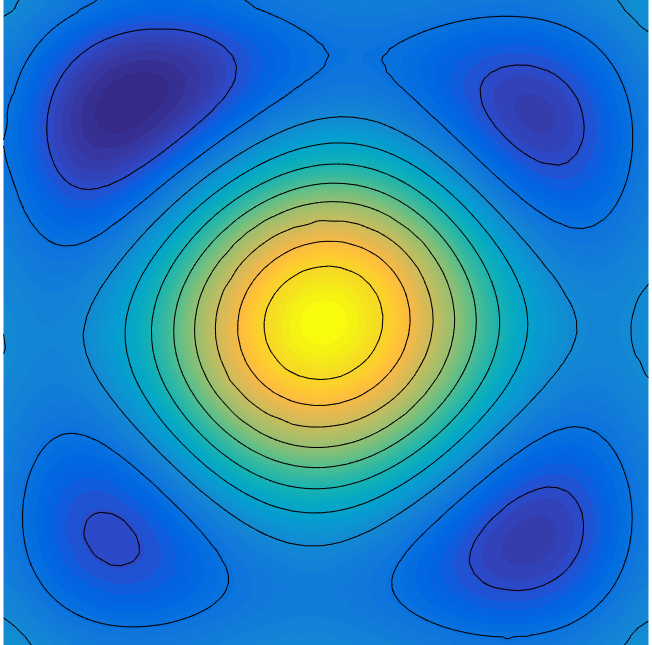}\\[0.15cm]

\rotatebox[origin=l]{90}{Scenario 1}& 
\includegraphics[width=0.115\textwidth]{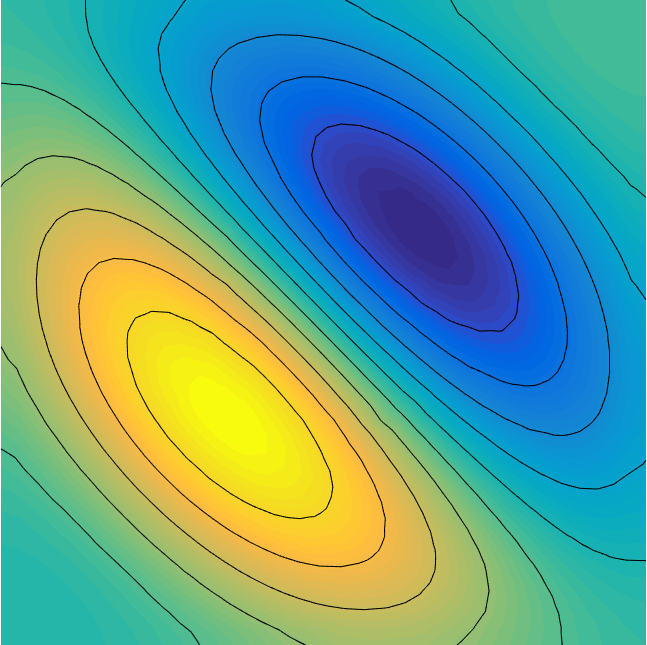} & 
\includegraphics[width=0.115\textwidth]{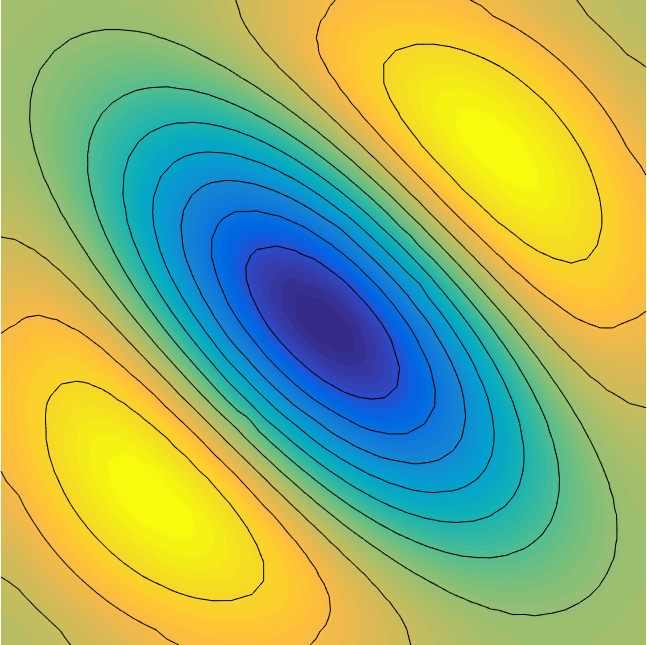} & 
\includegraphics[width=0.115\textwidth]{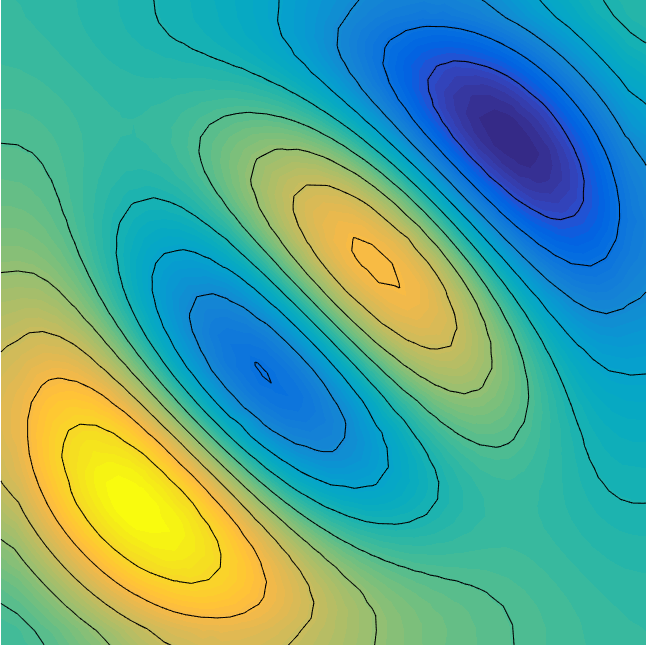} & 
\includegraphics[width=0.115\textwidth]{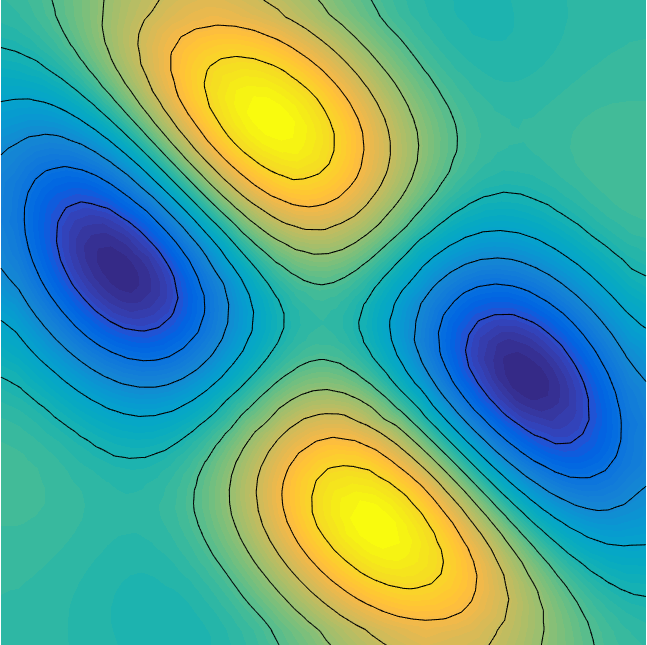} & 
\includegraphics[width=0.115\textwidth]{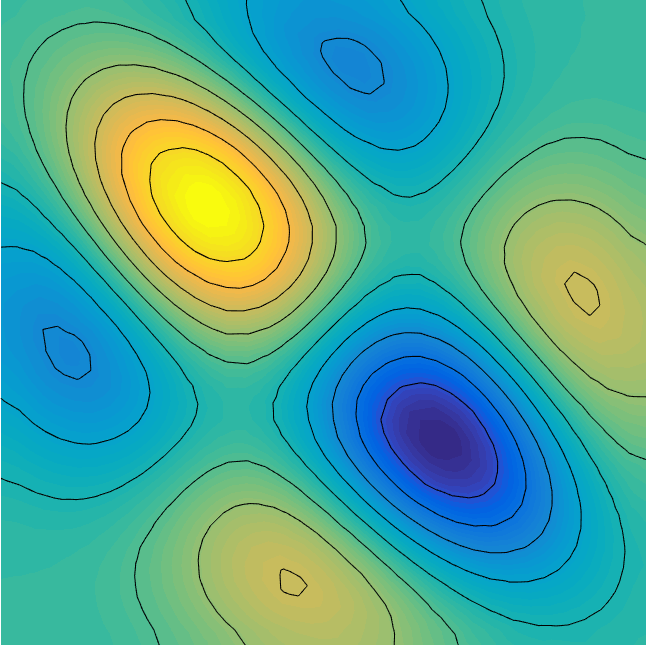} & 
\includegraphics[width=0.115\textwidth]{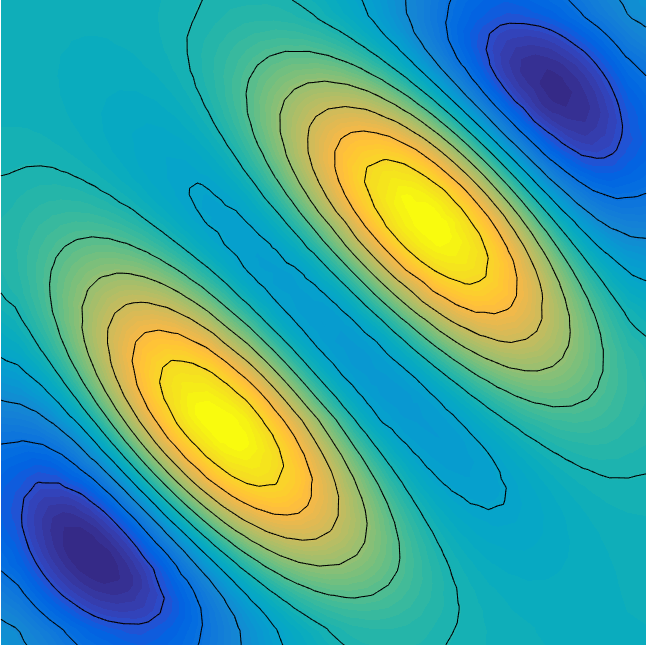}\\[0.15cm]

\rotatebox[origin=l]{90}{Scenario 2}& 
\includegraphics[width=0.115\textwidth]{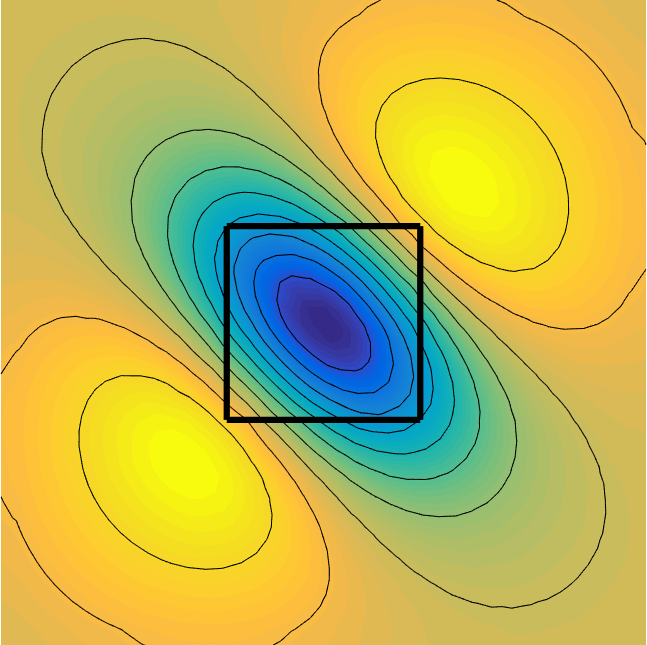} & 
\includegraphics[width=0.115\textwidth]{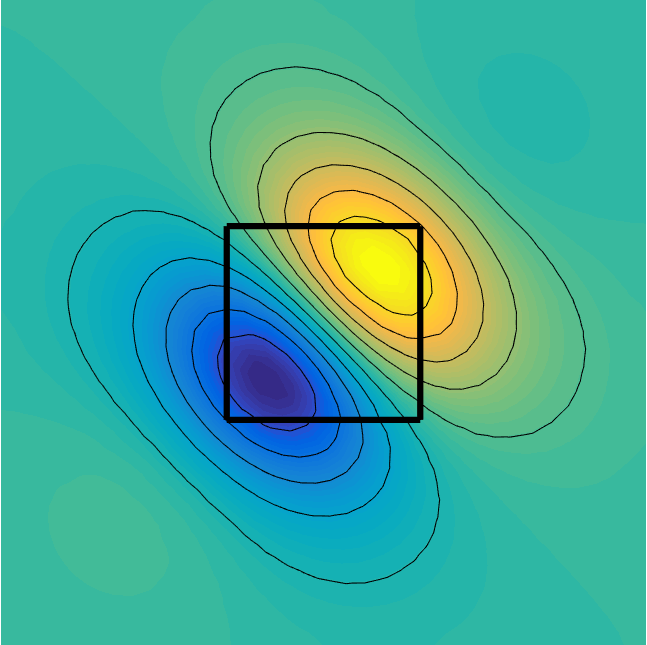} & 
\includegraphics[width=0.115\textwidth]{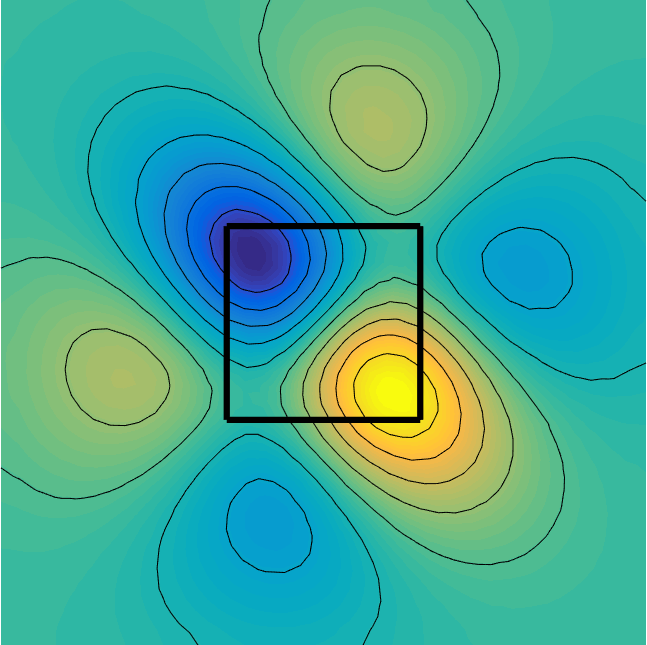} & 
\includegraphics[width=0.115\textwidth]{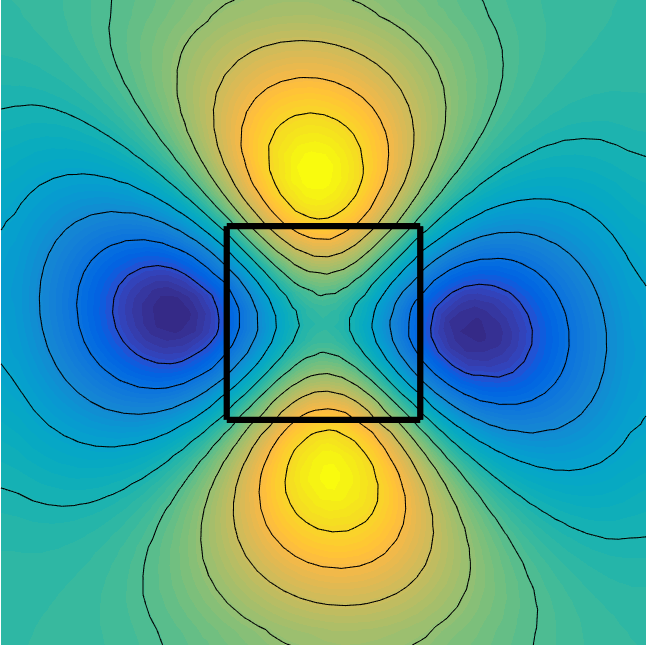} & 
\includegraphics[width=0.115\textwidth]{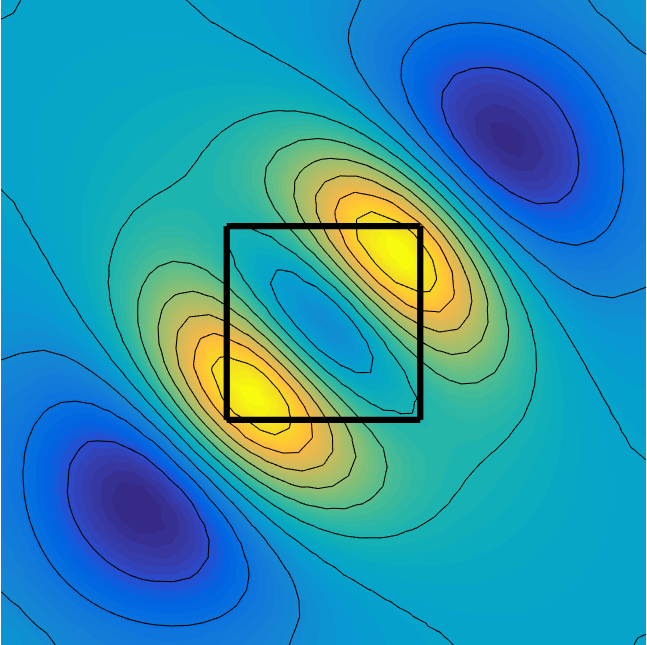} & 
\includegraphics[width=0.115\textwidth]{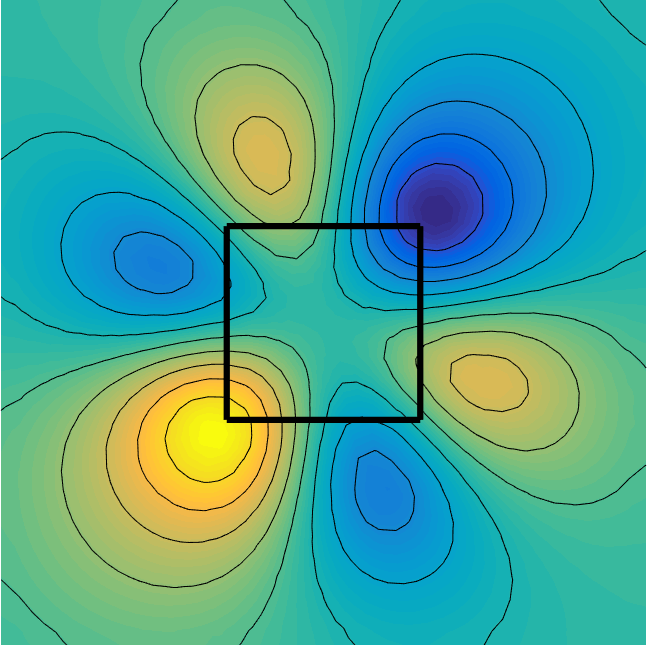}\\[0.15cm]

\rotatebox[origin=l]{90}{$\begin{array}{c} \text{Scenario 3} \\ \alpha=\beta=1 \end{array}$}& 
\includegraphics[width=0.115\textwidth]{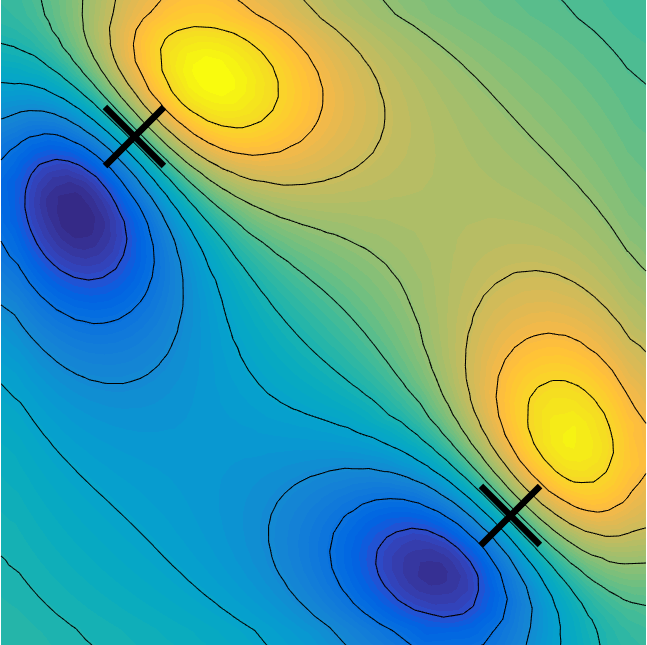} & 
\includegraphics[width=0.115\textwidth]{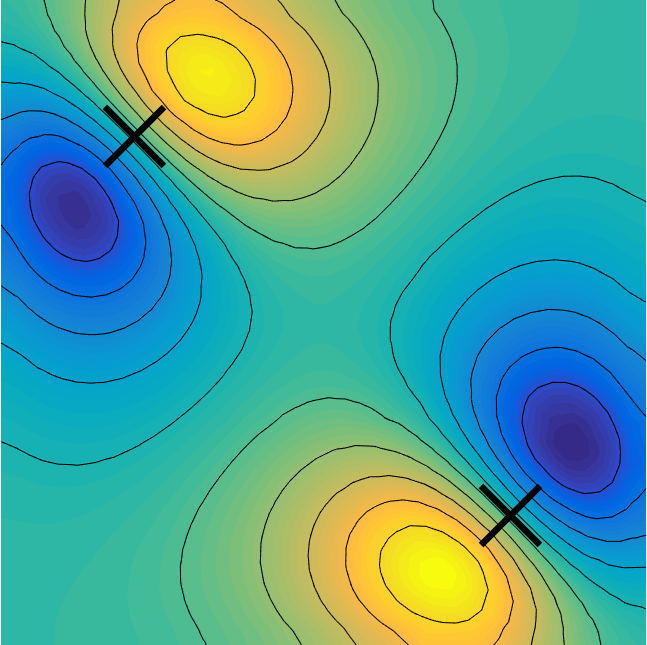} & 
\includegraphics[width=0.115\textwidth]{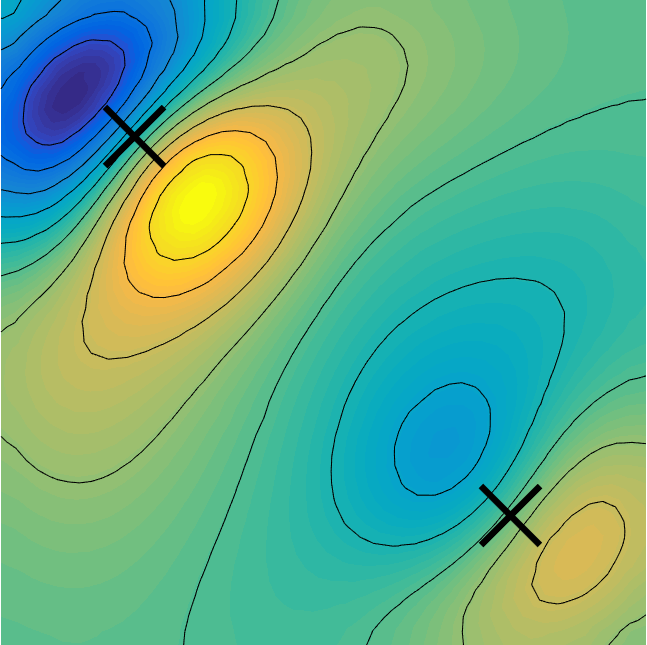} & 
\includegraphics[width=0.115\textwidth]{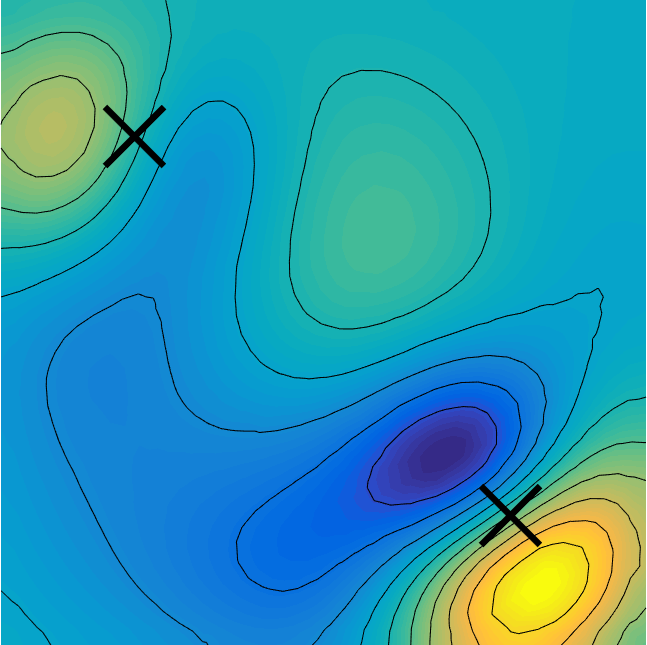} & 
\includegraphics[width=0.115\textwidth]{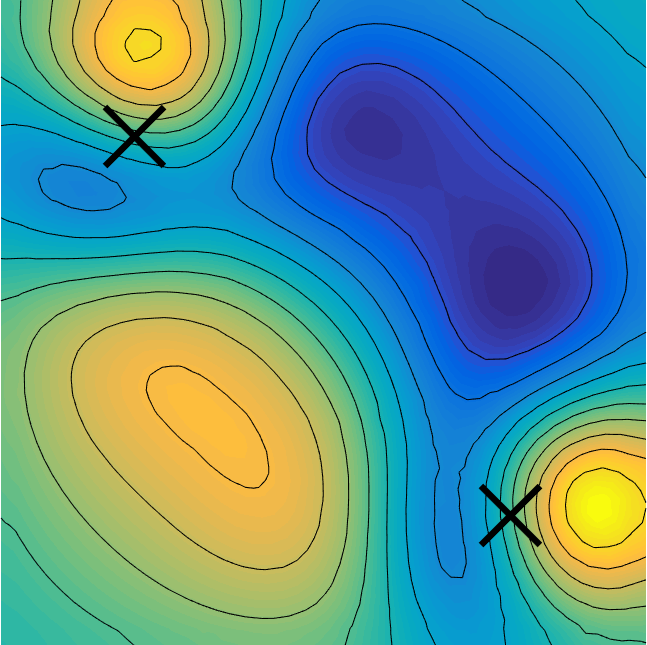} & 
\includegraphics[width=0.115\textwidth]{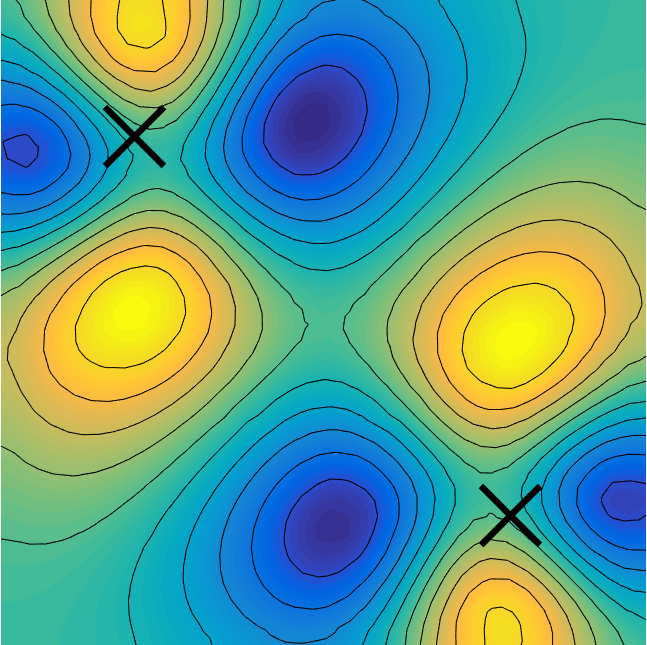} \\[0.15cm]

\rotatebox[origin=l]{90}{\hspace{-0.5cm}$\begin{array}{c} \text{Scenario 3} \\ \alpha=10,\beta=1 \end{array}$}& 
\includegraphics[width=0.115\textwidth]{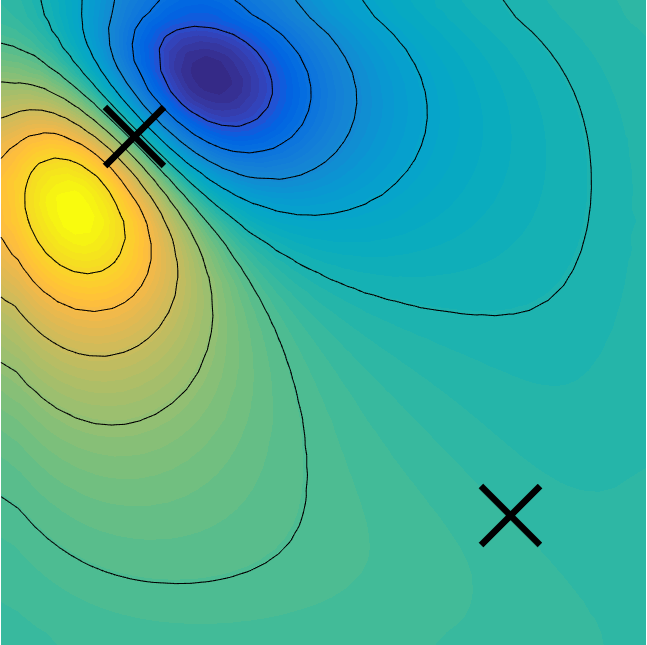} & 
\includegraphics[width=0.115\textwidth]{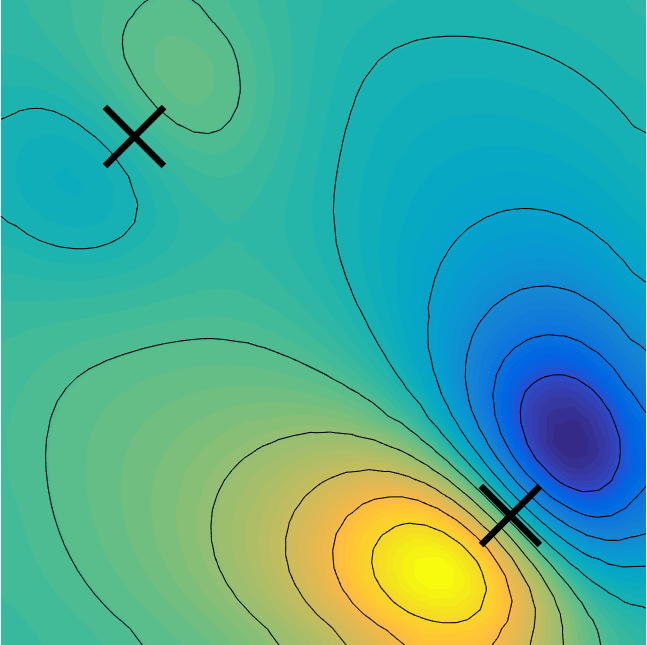} & 
\includegraphics[width=0.115\textwidth]{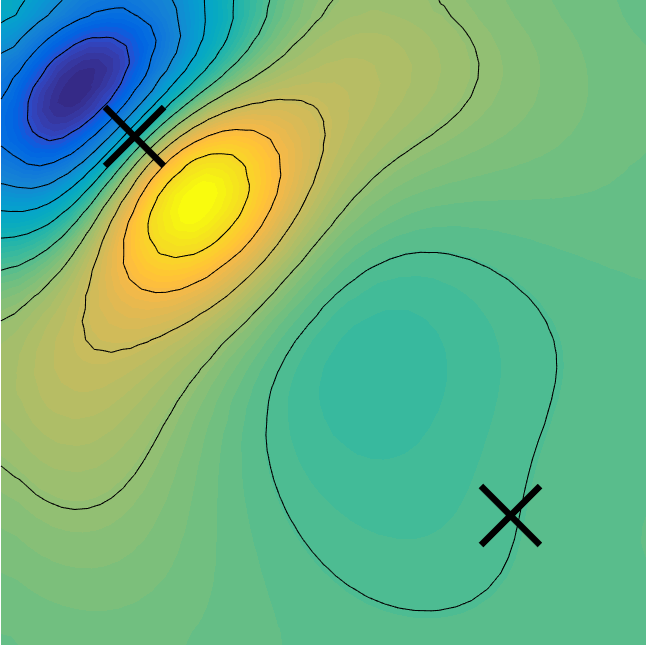} & 
\includegraphics[width=0.115\textwidth]{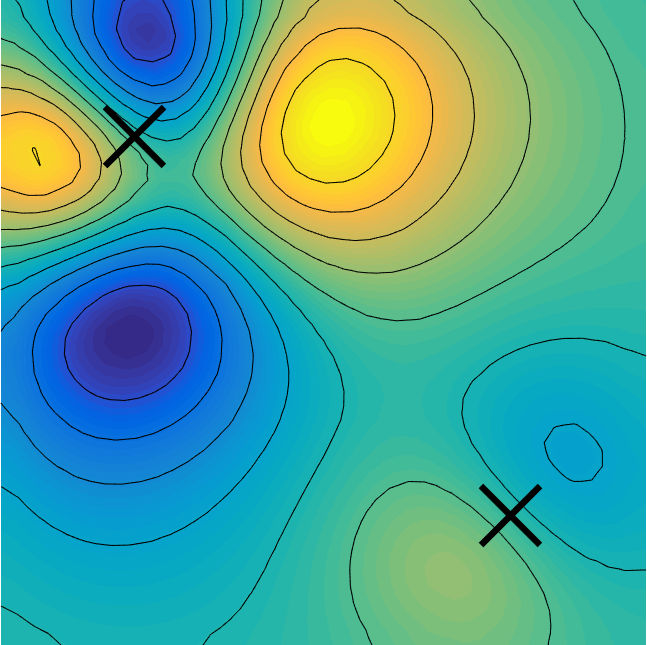} & 
\includegraphics[width=0.115\textwidth]{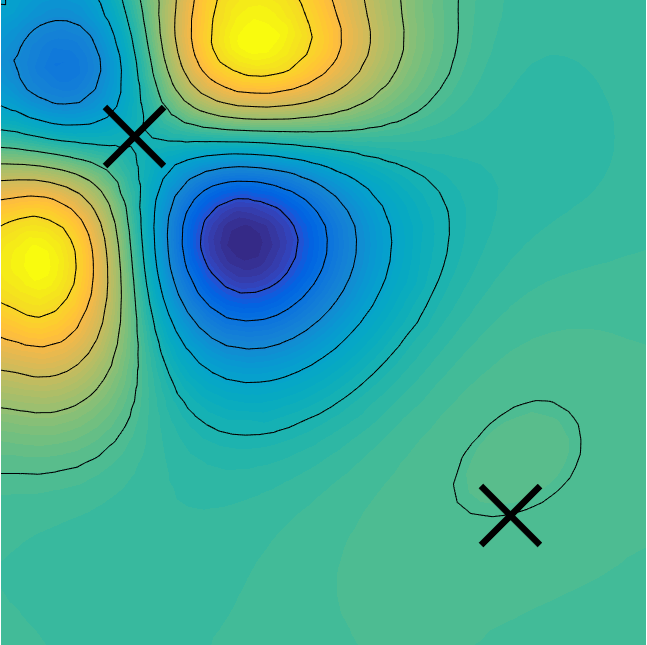} & 
\includegraphics[width=0.115\textwidth]{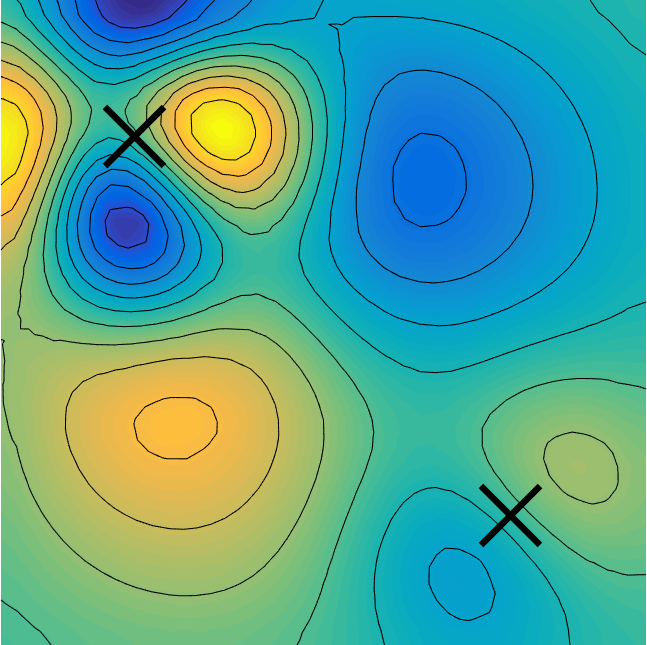}
\end{tabular}

\caption{Parameter modes: each figure represents the function $s\mapsto \sum_{i=1}^d v_i \, \mathbf{1}_i(s)$ for different $v\in\mathbb{R}^d$, where $\mathbf{1}_i$ is the indicator function of the $i$-th element of the mesh.
In the first row (K-L) $v$ is the $i$-th eigenvector of $\Sigma$, which corresponds to the Karhunen-Loève modes. In the four other rows, $v$ is the $i$-th generalized eigenvector of the matrix pair $(H,\Sigma^{-1})$, for different $H$ depending on the scenario.}

\label{fig:modes}
\end{figure}

For each scenario, an approximation $\widehat H$ of $H$ is computed with a large number of samples, $K=10^4$. This approximation is considered sufficiently accurate and will be used in place of $H$. Figure \ref{fig:modes} illustrates the leading generalized eigenvectors of the matrix pair $(H,\Sigma^{-1})$ as well as the leading eigenvectors of $\Sigma$, meaning the K-L modes; see Section \ref{sec:KL}. Since they do not depend upon $f$, the K-L modes do not have any particular relation to the elliptic PDE solution other than some symmetry properties related to the shape of the domain $\Omega$. In contrast, the modes associated with the three scenarios present specific features which depend on the function $f$. For example with scenario 2, we observe that the modes in the parameter space somehow represent more information \emph{local} to the region of interest $\Omega_s$.

The choice of the norm $\|\cdot\|_V$ also impacts the generalized eigenvectors of $(H,\Sigma^{-1})$ through the matrix $H$.
For instance with scenario 3 we have $R_V=\text{diag}(\alpha,\beta)$ allows us to write
$$
 H = \alpha \, H_1 + \beta \, H_2 \,,
 \quad\text{with}\quad 
 H_i=\int (\nabla f_i)^T(\nabla f_i)\mathrm{d}x \quad i=1,2.
$$
With the choice $\alpha=\beta=1$, the modes in Figure \ref{fig:modes} suggest that the two points of interest $s_a$ and $s_b$ are considered equally important, whereas the choice $\alpha=10 $ and $\beta=1$ leads to significantly more patterns around point $s_a$ (on the top-left of $\Omega$) than around point $s_b$.

\subsubsection{Approximating the conditional expectation and comparison with K-L}

\begin{figure}\centering
 \begin{subfigure}[b]{0.31\textwidth}\centering
  \includegraphics[width=0.95\textwidth]{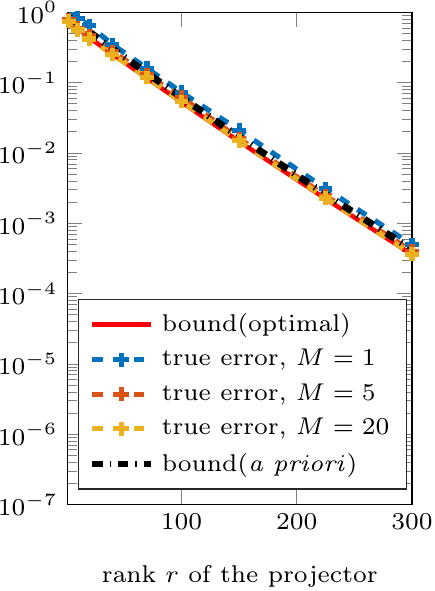}
  \caption{Scenario 1}
 \end{subfigure}
 \begin{subfigure}[b]{0.31\textwidth}\centering
  \includegraphics[width=0.95\textwidth]{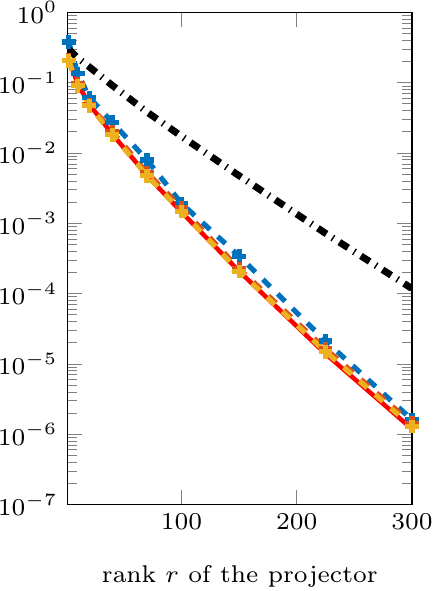}
  \caption{Scenario 2}
 \end{subfigure}
 \begin{subfigure}[b]{0.31\textwidth}\centering
  \includegraphics[width=0.95\textwidth]{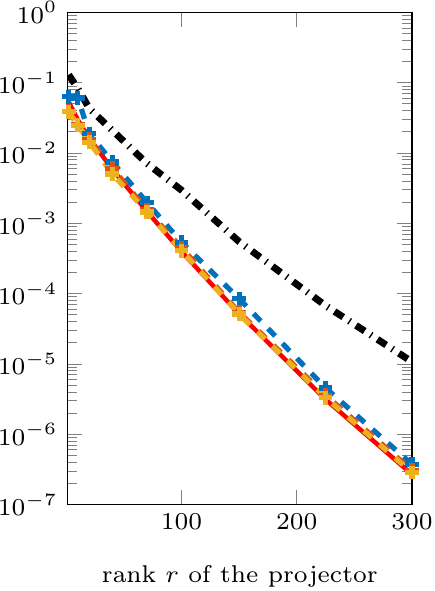}
  \caption{Scenario 3 ($\alpha=\beta=1$)}
 \end{subfigure}
 \caption{
 Error $\| f- \widehat F_r \|_\mathcal{H}$ as a function of the rank of $P_r$.
 The error $\| f- \widehat F_r \|_\mathcal{H}=\mathbb{E}(\|f(X)-\widehat F_r(X)\|_V^2)^{1/2}$, $X\sim\mu$, is estimated via Monte Carlo with 300 samples for $X$.
 The red (solid) and black (dash-dot) lines represent the upper bound $\trace( \Sigma(I_d-P_r^T) H (I_d-P_r) )^{1/2}$ with $P_r$ defined either as the minimizer of the upper bound (red lines) or as the projector onto the leading eigenspace of $\Sigma$ (black lines).
 }
 \label{fig:CondExp}
\end{figure}

Assume the matrix $H$ is known (again, a sufficiently accurate approximation $\widehat H$ with $K=10^4$ samples is used in place of $H$) and let $P_r$ be the rank-$r$ projector which minimizes $\trace( \Sigma(I_d-P_r^T) H (I_d-P_r) )$. We consider the approximation $\widehat F_r$ of the conditional expectation $\mathbb{E}_\mu(f|\sigma(P_r))$ given by \eqref{eq:defFhatr} with $\widehat P_r = P_r$. Figure \ref{fig:CondExp} shows the error $\| f- \widehat F_r \|_\mathcal{H}$ as a function of the rank $r$ of the projector. For each scenario, one realization of $\widehat F_r$ is computed with either $M=1$, $M=5$, or $M=20$ samples. We first note that, since we do not exactly compute the conditional expectation, the errors (dotted curves) are sometimes above the upper bound (solid red curves). In this inexact setting, $\trace( \Sigma(I_d-P_r^T) H (I_d-P_r) )^{1/2}$ is no longer a certified upper bound for the error. However we observe it can still be used as a good error indicator.

The three scenarios do not have the same convergence rate with $r$: the first scenario has the slowest and the third the fastest. Even though they are different post-solution treatments of the same solution map $x\mapsto u(x)$, the functions $f$ associated with each scenario do not have the same complexity in terms of intrinsic dimension. 

Interestingly, increasing $M$ does not lead to significant improvements of the approximation. This phenomenon can be explained by the following relation,
\begin{align*}
 \mathbb{E} \big(  \| f - \widehat F_r \|_{\mathcal{H}}^2  \big)
 &\overset{\eqref{eq:CondExp_VarForm}}{=}\mathbb{E} \big(  \| f -\mathbb{E}_\mu(f|\sigma(P_r)) \|_{\mathcal{H}}^2 + \| \widehat F_r - \mathbb{E}_\mu(f|\sigma(P_r)) \|_{\mathcal{H}}^2  \big) \\
 &\overset{\eqref{eq:defFhatr}}{=} \| f -\mathbb{E}_\mu(f|\sigma(P_r)) \|_{\mathcal{H}}^2  + \frac{1}{M} \| f - \mathbb{E}_\mu(f|\sigma(P_r)) \|_{\mathcal{H}}^2  \\
 &= \Big(1+\frac{1}{M}\Big) \| f - \mathbb{E}_\mu(f|\sigma(P_r)) \|_{\mathcal{H}}^2 \,,
\end{align*}
where the expectation is taken over the samples $Y_1,\hdots,Y_M$ (the projector $P_r$ being fixed here). This result shows that even with small $M$, one can still hope to obtain a good approximation $\widehat F_r$ of $f$ provided $P_r$ is chosen such that $\| f - \mathbb{E}_\mu(f|\sigma(P_r)) \|_{\mathcal{H}}$ is sufficiently small. In other words a crude approximation of the conditional expectation yields at most a factor of two (when $M=1$) in the expected error squared, so that it remains of the same order of magnitude as $\| f - \mathbb{E}_\mu(f|\sigma(P_r)) \|_{\mathcal{H}}^2$; see also Theorem 3.2 from \cite{Constantine2014a}.

\begin{figure}\centering
 \scriptsize
 \setlength\figwidth{0.3\textwidth}
 \setlength\figheight{0.3\textwidth}
 \begin{subfigure}[b]{0.4\textwidth}\centering
  \includegraphics[width=0.95\textwidth]{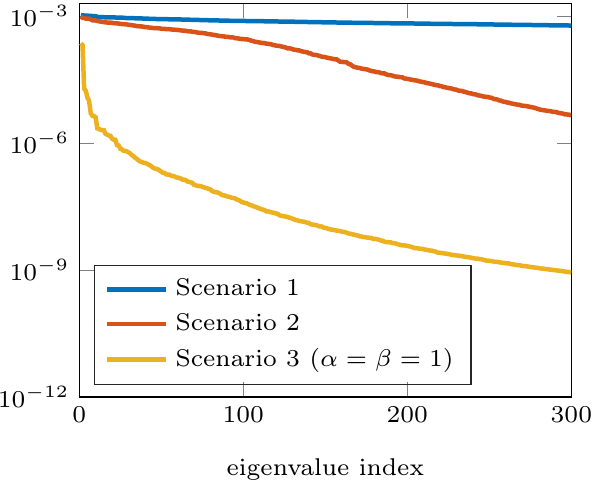}
 \end{subfigure}
 \begin{subfigure}[b]{0.4\textwidth}\centering
  \includegraphics[width=0.95\textwidth]{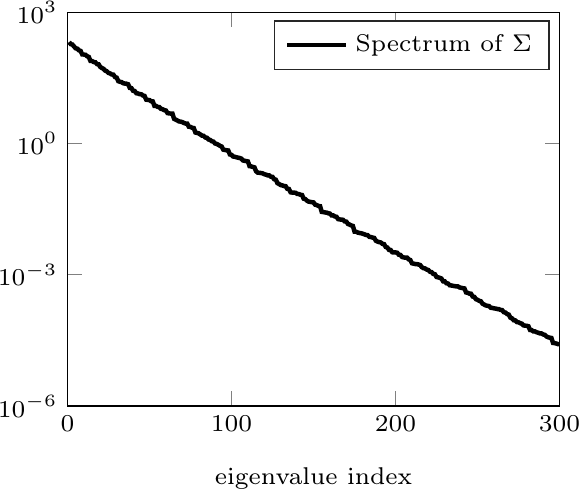}
 \end{subfigure}
 \caption{Spectrum of $H$ for the three scenarios (left) and spectrum of $\Sigma$ (right).}
 \label{fig:spectra}
\end{figure}

We now compare with the truncated Karhunen-Loève decomposition, for which $P_r$ is defined as the rank-$r$ orthogonal projector onto the leading eigenspace of the covariance matrix $\Sigma$. The black dash-dotted curves in Figure \ref{fig:CondExp} represent the upper bound $\trace( \Sigma(I_d-P_r^T) H (I_d-P_r) )^{1/2}$ for this choice of $P_r$, as a function of $r$. (The true error $\| f - \mathbb{E}_\mu(f|\sigma(P_r)) \|_{\mathcal{H}}$ is substantively the same as its upper bound, so we decided not to plot it.) It is interesting to see that in the first scenario, the K-L projector is essentially as effective as the projector obtained by minimizing the upper bound. As shown in Figure \ref{fig:spectra}, the spectrum of $H$ is flat, which means that $H$ is close to a rescaled identity matrix. Then, minimizing $\trace( \Sigma(I_d-P_r^T) H (I_d-P_r) )$ is nearly the same as minimizing $\trace( (I_d-P_r)\Sigma(I_d-P_r^T) )=\mathbb{E}(\|X-P_rX\|_2^2)$, where $X\sim\mathcal{N}(0,\Sigma)$, and yields the same projector as the truncated K-L method; see \eqref{eq:KLmotivation}.
However, this reasoning does not apply to scenarios 2 and 3, where the spectrum of $H$ decays rapidly. For these scenarios we observe in Figure \ref{fig:CondExp} that the new method outperforms the truncated K-L method. For instance, in scenario 2 the new method reaches an error of $10^{-4}$ with only $r=150$ whereas the truncated K-L method requires $r=300$.

\subsubsection{Quality of the projector}

\begin{figure}\centering
 \scriptsize
 \setlength\figwidth{0.235\textwidth}
 \setlength\figheight{0.25\textheight}
 \begin{subfigure}[b]{0.31\textwidth}\centering
  \includegraphics[width=0.95\textwidth]{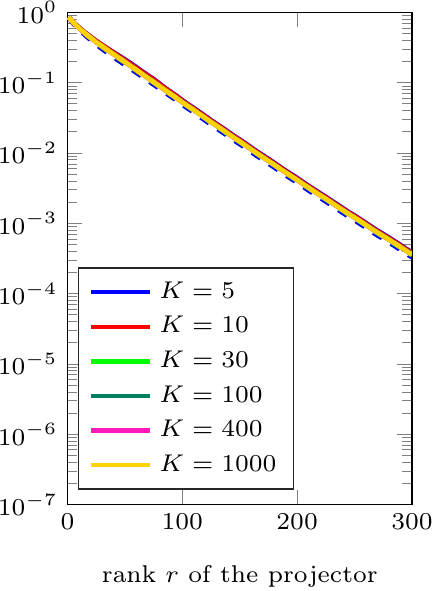}
  \caption{Scenario 1}
 \end{subfigure}
 \begin{subfigure}[b]{0.31\textwidth}\centering
  \includegraphics[width=0.95\textwidth]{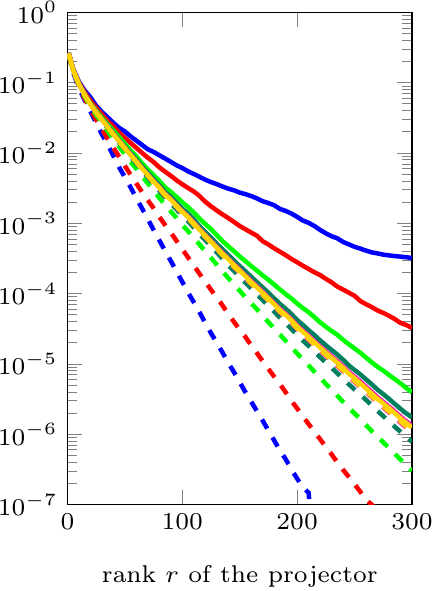}
  \caption{Scenario 2}
 \end{subfigure}
 \begin{subfigure}[b]{0.31\textwidth}\centering
  \includegraphics[width=0.95\textwidth]{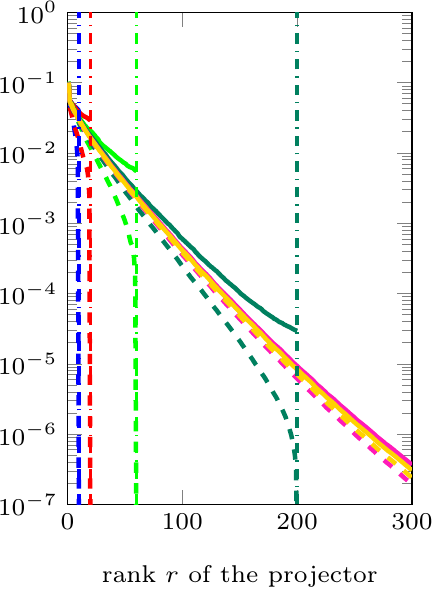}
  \caption{Scenario 3 ($\alpha=\beta=1$)}
  \label{fig:projector3}
 \end{subfigure}
 \caption{
 Error bound $\trace( \Sigma(I_d-\widehat P_r^T) H (I_d-\widehat P_r) )^{1/2}$ (solid curves) and approximate error bound $\trace( \Sigma(I_d-\widehat P_r^T) \widehat H (I_d-\widehat P_r) )^{1/2}$ (dotted curves) as a function of the rank of $\widehat P_r$. For each scenario, the curves correspond to one realization of $\widehat H$ and $\widehat P_r$ defined by \eqref{eq:defHhat} and \eqref{eq:defPrHat} for different values of $K$. In Figure \ref{fig:projector3}, the vertical lines correspond to $r=2K$.
 }
 \label{fig:projector}
\end{figure}

In this section we assess the quality of a projector $\widehat P_r$ defined by \eqref{eq:defPrHat}, where $\widehat H$ is the $K$-sample Monte Carlo approximation of $H$ given by \eqref{eq:defHhat}. In the present context, an \emph{optimal} projector would be a minimizer of $P_r\mapsto\trace( \Sigma(I_d-P_r^T) H (I_d-P_r) )^{1/2}$ so that the only relevant criteria for the quality of $\widehat P_r$ is how close $\trace( \Sigma(I_d-\widehat P_r^T) H (I_d-\widehat P_r) )^{1/2}$ is to the minimum of the upper bound. Figure \ref{fig:projector} contains two sets of curves: the solid curves represent the error bound $\trace( \Sigma(I_d-\widehat P_r^T) H (I_d-\widehat P_r) )^{1/2}$ as a function of the rank of $\widehat P_r$, whereas the dotted curves correspond to the \emph{approximate} error bound $\trace( \Sigma(I_d-\widehat P_r^T) \widehat H (I_d-\widehat P_r) )^{1/2}$. This approximate error bound is the quantity we would use in place of the error bound when the matrix $H$ is not known. For each scenario we observe that for small $K$, the approximate error bound underestimates the true error bound. This means that $\trace( \Sigma(I_d-\widehat P_r^T) \widehat H (I_d-\widehat P_r) )^{1/2}$ can be used as an error estimator only if $K$ is sufficient large.

Observe in Figure \ref{fig:projector} that scenarios 1 and 2 need fewer samples to obtain a good projector (say around $K=30$ samples) compared to the last scenario (at least $K=400$ samples). To understand this result, let us note that if $r$ is larger than the rank of $\widehat H$, the projector $\widehat P_r$ is not uniquely determined: any $\widehat P_r$ such that $\text{Im}(\widehat H) \subset \text{Im}(\widehat P_r)$ is a solution to \eqref{eq:defPrHat}. Therefore the rank of $\widehat P_r$ should not exceed that of $\widehat H$ which, thanks to \eqref{eq:defHhat}, satisfies the following relation
$$
 \text{rank}(\widehat H) \leq K\, \text{rank}\Big( \big(\nabla f(X) \big)^TR_V \big( \nabla f(X) \big) \Big) \leq K\, \text{dim}(V) \,.
$$
With scenario 3 we have $\text{dim}(V)=2$ so that the rank of $\widehat P_r$ should not exceed $2K$. This limitation is represented by the vertical lines on Figure \ref{fig:projector3}. With scenarios 1 and 2 we have $\text{dim}(V)=1691$ and  $\text{dim}(V)=168$ so that this limit is not attained within the range of the plots. The conclusion is that when the dimension of $V$ is large, one needs fewer samples from $\nabla f(X)$ to obtain a suitable projector, because each sample is a matrix with potentially a large rank.

\section{Conclusions}

We have addressed the problem of approximating multivariate functions taking values in a vector space. We approximate such functions by means of \emph{ridge functions} that depend on a number of linear combinations of the input parameters that is smaller than the original dimension. Rather than seeking an optimal approximation, we build a \emph{controlled} approximation: we develop an upper bound on the approximation error and minimize this upper bound.

Our analytical and numerical examples demonstrate good performance of the method, and also illustrate conditions under which it might not work well. For example, we show cases where minimizing the upper bound leads to an optimal approximation, and contrasting cases where the error bound is not tight. Numerical demonstrations on an elliptic PDE also illustrate various computational issues: sampling to compute both the projector (yielding the important directions) and the conditional expectation (yielding the ridge profile). 

Future work may explore several natural extensions of the proposed methodology. First is the extension to non-Gaussian input measures, e.g., uniform measure on bounded domains in $\mathbb{R}^d$. Second is the extension to infinite-dimensional input spaces: for example, letting the domain of $f$ be a function space endowed with Gaussian or Besov measure. Finally, it may be possible to develop sharper error bounds based on higher-order derivatives, e.g., Hessians of $f$. For the last two points, we may be able to use recent results on higher-order Poincaré inequalities \cite{nourdin2009second}. 

\appendix
\section{Proofs}\label{sec:proofs}

\subsection{Proof of Proposition \ref{prop:ImportanceKernel}}\label{proof:prop:ImportanceKernel}
 Let $h \in \mathcal{H}_{P_r}$. By \eqref{eq:Doob} we can write $h=g\circ P_r$ for some Borel function $g$. Since $\text{Ker}(Q_r)=\text{Ker}(P_r)$ we have $P_rx = P_rQ_rx =0$ for all $x\in\text{Ker}(Q_r)$. Also for any $x\in\text{Im}(Q_r)$ we have $Q_rx=x$ and then $P_rx = P_rQ_rx $. Thus $P_rx = P_rQ_rx$ holds for any $x\in\mathbb{R}^d = \text{Ker}(Q_r)\oplus \text{Im}(Q_r)$ so that $P_r = P_r Q_r$. Then $h=g\circ P_r = (g\circ P_r )\circ Q_r$ which shows that $h\in \mathcal{H}_{Q_r}$. 
 Then the inclusion $\mathcal{H}_{P_r} \subset \mathcal{H}_{Q_r}$ holds. By symmetry of the role of $P_r$ and $Q_r$ we obtain the result.

\subsection{Proof of Proposition \ref{prop:ExplicitCondExp}}\label{proof:prop:ExplicitCondExp}

 Let $F:x\mapsto \int_{\mathbb{R}^d} f(P_r x + (I_d-P_r)y) \mu(\mathrm{d}y)$ and $h\in\mathcal{H}_{P_r}$. By \eqref{eq:Doob}, $h$ can be written as $g\circ P_r$ for some Borel function $g$ so that $h(x)=h(P_rx+(I_d-P_r)y)$ for all $x,y\in\mathbb{R}^d$. We can write
 \begin{align*}
  \int_{\mathbb{R}^d} & ( F(x),h(x) )_V \mathrm{d}\mu(x)
  =\int_{\mathbb{R}^d} \Big( \int_{\mathbb{R}^d} f(P_r x + (I_d-P_r)y) \mu(\mathrm{d}y),h(x) \Big)_V \mathrm{d}\mu(x) \\
  &=\int_{\mathbb{R}^d} \int_{\mathbb{R}^d} \Big( f(P_r x + (I_d-P_r)y) ,h(P_r x + (I_d-P_r)y) \Big)_V \mu(\mathrm{d}y)\mathrm{d}\mu(x) \\
  &=\mathbb{E}\big( ( f(Z),h(Z) )_V \big)  \, ,
 \end{align*}
 where the expectation is taken over the random vector $Z=P_rX+(I_d-P_r)Y$, where $X$ and $Y$ are two independent random vectors distributed as $\mu=\mathcal{N}(m,\Sigma)$. If $Z\sim\mu$ then the previous relation yields \eqref{eq:CondExp_VarForm} for any $h\in\mathcal{H}_{P_r}$, which would conclude the proof.
 
 It remains to show that $Z\sim\mu$. Note that $Z$ is Gaussian with mean $m$ and covariance 
 $$
  \text{Cov}(Z) = P_r\Sigma P_r^T +(I_d-P_r)\Sigma (I_d-P_r^T) = \Sigma - P_r\Sigma - \Sigma P_r^T + 2 P_r \Sigma P_r^T \, .
 $$
 Then $Z\sim\mu$ if and only if $P_r\Sigma + \Sigma P_r^T = 2 P_r \Sigma P_r^T$.
 Since $P_r$ is $\Sigma^{-1}$-orthogonal, relation \eqref{eq:orthogonalityOfPr} holds for any $x\in\mathbb{R}^d$ which is equivalent to
 $$
  P_r^T\Sigma^{-1} + \Sigma^{-1}P_r = 2 P_r^T\Sigma^{-1}P_r \,.
 $$
 Multiplying by $P_r^T$ to the left (resp.\ by $P_r$ to the right) we get $P_r^T\Sigma^{-1} = P_r^T\Sigma^{-1}P_r$ (resp.\ $\Sigma^{-1}P_r = P_r^T\Sigma^{-1}P_r$) so that the relation $P_r^T\Sigma^{-1} = \Sigma^{-1}P_r$ holds and yields $\Sigma P_r^T = P_r\Sigma = P_r\Sigma P_r^T$. Therefore we have
 $$
  P_r\Sigma + \Sigma P_r^T = 2 P_r \Sigma P_r^T\, ,
 $$
 which concludes the proof.

\subsection{Proof of Proposition \ref{prop:subspin}}\label{proof:prop:subspin}

 First we assume that $P_r$ is a $\Sigma^{-1}$-orthogonal projector. 
 Let $h:\mathbb{R}^d\rightarrow\mathbb{R}$ be a continuously differentiable function and define $g:x\mapsto h( P_r y + (I_d-P_r)x )$ for some $y\in\mathbb{R}^d$. For any $x\in\mathbb{R}^d$ we have $\nabla g( x ) = (I_d-P_r)^T\nabla h( P_r y + (I_d-P_r)x )$. By Proposition \ref{prop:ExplicitCondExp} we have
 $$
  \mathbb{E}_\mu(g) = \int_{\mathbb{R}^d} h( P_r y + (I_d-P_r)x' )\mu(\mathrm{d}x') = \mathbb{E}_\mu(h|\sigma(P_r))(y) \, .
 $$
 Notice that we can write $\mathbb{E}_\mu(h|\sigma(P_r))(y) = \mathbb{E}_\mu(h|\sigma(P_r))(P_r y + (I_d-P_r)x)$. Then the Poincaré inequality \eqref{eq:Poincare} applied with the function $g$ yields
 \begin{align*}
  \int_{\mathbb{R}^d} \big(  h(P_r y  + (I_d-P_r)x) - &\mathbb{E}_\mu(h|\sigma(P_r))(P_r y + (I_d-P_r)x) \big)^2 \mathrm{d}\mu(x) \\
  &\leq \int_{\mathbb{R}^d} \| (I_d-P_r)^T\nabla h( P_r y + (I_d-P_r)x ) \|_{\Sigma}^2 ~ \mathrm{d}\mu(x) \,.
 \end{align*}
 Recall that, since $P_r$ is $\Sigma^{-1}$-orthogonal, we have $P_r Y + (I_d-P_r)X \sim\mu $ whenever $X\sim\mu$ and $Y\sim\mu$ are independent; see the proof of Proposition \ref{prop:ExplicitCondExp}. Thus, replacing $y$ by $Y$ in the previous inequality and taking the expectation over $Y$ yields \eqref{eq:SubspacePoincare}.

 It remains to show that \eqref{eq:SubspacePoincare} also holds for projectors that are not $\Sigma^{-1}$-orthogonal. Thus let $P_r$ be any projector and define $Q_r$ as the (unique) $\Sigma^{-1}$-orthogonal projector such that $\text{Ker}(Q_r) = \text{Ker}(P_r)$. 
 Following the proof of Proposition \ref{prop:ExplicitCondExp}, we have that $Q_r$ satisfies $Q_r\Sigma + \Sigma Q_r^T = 2 Q_r \Sigma Q_r^T$ which is equivalent to saying that the relation $\|x\|_{\Sigma}^2 = \|Q_r^Tx\|_{\Sigma}^2 + \|(I_d-Q_r^T)x\|_{\Sigma}^2 $ holds for any $x\in\mathbb{R}^d$. Then $\|x\|_{\Sigma}^2 \geq \|(I_d-Q_r^T)x\|_{\Sigma}^2$ for any $x\in\mathbb{R}^d$. Replacing $x$ by $(I_d-P_r^T)x$ we get
 \begin{align}
  \|(I_d-P_r^T)x\|_{\Sigma}^2 
  &\geq \|(I_d-Q_r^T)(I_d-P_r^T)x\|_{\Sigma}^2 \nonumber\\
  &= \|(I_d-Q_r^T-P_r^T+ Q_r^TP_r^T)x\|_{\Sigma}^2 \nonumber\\
  &= \|(I_d-Q_r^T)x\|_{\Sigma}^2 \,.  \label{eq:tmp326872}
 \end{align}
 For the last equality we used relation $P_r=P_rQ_r$, which holds true since $\text{Ker}(P_r)=\text{Ker}(Q_r)$.
 Finally, Proposition \ref{prop:ImportanceKernel} allows writing $\mathbb{E}_\mu(h|\sigma(P_r)) = \mathbb{E}_\mu(h|\sigma(Q_r))$ so that
 \begin{align*}
  \int ( h - \mathbb{E}_\mu(h & |\sigma(P_r)) )^2 \, \mathrm{d}\mu
  = \int ( h - \mathbb{E}_\mu(h|\sigma(Q_r)) )^2 \, \mathrm{d}\mu \\
  &\overset{\eqref{eq:SubspacePoincare}}{\leq} \int \| (I_d-Q_r^T)\nabla h \|_{\Sigma}^2 \, \mathrm{d}\mu 
  \overset{\eqref{eq:tmp326872}}{\leq} \int \| (I_d-P_r^T)\nabla h \|_{\Sigma}^2 \, \mathrm{d}\mu \,,
 \end{align*}
 which shows that \eqref{eq:SubspacePoincare} holds for any projector $P_r$.

\subsection{Proof of Proposition \ref{prop:BoundVVfunction}}\label{proof:prop:BoundVVfunction}

 Denote by $(w_i,\alpha_i)\in\mathbb{R}^n\times \mathbb{R}_{\geq0}$ the $i$-th eigenpair of the matrix $R_V$ so that $R_V=\sum_{i=1}^n \alpha_i w_iw_i^T $ and $\|y\|_V^2 = \sum_{i=1}^n \alpha_i (w_i^T y)^2$ for any $y\in V$.
 The function $f$ can be represented as $x\mapsto \sum_{i=1}^n f_i(x) w_i $ where $f_i : x \mapsto w_i^T f(x)$.
 The linearity of the conditional expectation permits to write
 \begin{align}
  \| f - \mathbb{E}_\mu(f|\sigma(P_r)) \|_{\mathcal{H}}^2 
  &= \sum_{i=1}^n \alpha_i \int ( f_i - \mathbb{E}_\mu(f_i|\sigma(P_r)) )^2 \, \mathrm{d}\mu . \label{eq:tmp5328}
 \end{align}
 Because $f$ is continuously differentiable, the coordinate $f_i$ are continuously differentiable as well. Then the subspace Poincaré inequality \eqref{eq:SubspacePoincare} yields
 \begin{align*}
  \int ( f_i - \mathbb{E}_\mu(f_i|\sigma(P_r)) )^2 \mathrm{d}\mu
  &\leq \int \| (I_d-P_r^T)\nabla f_i \|_{\Sigma}^2 ~ \mathrm{d}\mu \\
  &= \int \trace\big( \Sigma(I_d-P_r^T)(\nabla f_i)(\nabla f_i)^T (I_d-P_r) \big)  \mathrm{d}\mu \\
  &= \trace\big( \Sigma(I_d-P_r^T) \Big(\int  (\nabla f_i) (\nabla f_i)^T  \mathrm{d}\mu \Big) (I_d-P_r) \big).
 \end{align*}
 By definition of the Jacobian matrix \eqref{eq:defNablaF} we have $\nabla f_i(x) = \nabla f(x)^T w_i$. Then, together with \eqref{eq:tmp5328}, the above relation yields 
 \begin{align*}
  \| f & - \mathbb{E}_\mu(f|\sigma(P_r))  \|_{\mathcal{H}}^2 \\
  &\leq  \sum_{i=1}^n \alpha_i   \trace\big( \Sigma(I_d-P_r^T) \Big( \int  \nabla f(x)^T w_i w_i^T\nabla f(x)  \mathrm{d}\mu \Big) (I_d-P_r) \big) \\
  &= \trace\big( \Sigma(I_d-P_r^T) \Big( \int  \nabla f(x)^T \Big(\sum_{i=1}^n \alpha_i  w_i w_i^T\Big) \nabla f(x)  \mathrm{d}\mu \Big) (I_d-P_r) \big) \\
  &= \trace\big( \Sigma(I_d-P_r^T) H (I_d-P_r) \big),
 \end{align*}
 where for the last equality we used the relation $R_V=\sum_{i=1}^n \alpha_i w_iw_i^T $ and the definition \eqref{eq:defH} of $H$.
 This concludes the proof.

\subsection{Proof of Proposition \ref{prop:OptimalPr}} \label{proof:prop:OptimalPr}

 Let $H^{1/2}$ and $\Sigma^{1/2}$ be symmetric square roots of $H$ and $\Sigma$ respectively. For any projector $P_r$ we have
 $$
  \trace( \Sigma(I_d-P_r^T) H (I_d-P_r) )  = \| H^{1/2}(I_d-P_r)\Sigma^{1/2} \|_F^2 = \|A-X_r\|_F^2 \,,
 $$
 where $A=H^{1/2}\Sigma^{1/2}$ and $X_r=H^{1/2}P_r\Sigma^{1/2}$ and where $\|\cdot\|_F=\sqrt{\trace(\cdot)^T(\cdot)}$ denotes the Frobenius norm. Consider the singular value decomposition of $A=UDV^T$ where $U,V\in\mathbb{R}^{d\times d}$ are two orthogonal matrices and $D=\text{diag}(a_1,\hdots,a_d)$ with $a_1\geq a_2\geq\hdots\geq0$. The Eckart-Young theorem states that (i) the matrix $A_r = UD_rV^T$, with $D_r=\text{diag}(a_1,\hdots,a_r,0,\hdots,0)$, is a minimizer of $\|A-\widetilde A_r\|_F^2$ over all matrices $\widetilde A_r$ with $\text{rank}(\widetilde A_r)\leq r$ and (ii) that $\|A-A_r\|_F^2=a_{r+1}^2+\hdots+a_d^2$. 
 We now show that $A_r$ can be written as $X_r = H^{1/2}P_r\Sigma^{1/2}$ for some rank-$r$ projector $P_r$. Let $V_r\in\mathbb{R}^{d\times r}$ be the matrix containing the $r$ first columns of $V$ and let $P_r = \Sigma^{1/2} V_rV_r^T \Sigma^{-1/2}$. Since $V_r^TV_r=I_r$ we have $P_r^2=P_r$ so that $P_r$ is a rank-$r$ projector. Also we have $X_r = H^{1/2}P_r\Sigma^{1/2} = AV_rV_r^T = A_r$. Then $\|A-X_r\|^2=\|A-A_r\|^2 \leq \|A-\widetilde A_r\|^2$ holds for any rank-$r$ matrix $\widetilde A_r$, in particular for the ones of the form of $\widetilde A_r = H^{1/2}\widetilde P_r\Sigma^{1/2}$ for any rank-$r$ projector $\widetilde P_r$. This shows that the minimum in \eqref{eq:OptimalPr_residual} is reached by $P_r= \Sigma^{1/2} V_rV_r^T \Sigma^{-1/2}$. Furthermore it is easy to check that $P_r^T \Sigma^{-1} + \Sigma^{-1}P_r = 2P_r^T\Sigma^{-1}P_r$ holds so that, as we saw in the proof of Proposition \ref{prop:ExplicitCondExp}, $P_r$ is $\Sigma^{-1}$-orthogonal.
 
 It remains to show that $P_r$ can be written as in \eqref{eq:OptimalPr}.
 Notice that $A^TA=\Sigma^{1/2}H\Sigma^{1/2}=VD^2V^T$ holds and yields $H\Sigma^{1/2}V=\Sigma^{-1/2}VD^2$. 
 Denoting by $v_i$ the $i$-th column of $\Sigma^{1/2}V$ (which is such that $\|v\|_{\Sigma^{-1}}^2=1$), the latter relation yields $Hv_i=a_i^2\Sigma^{-1}v_i$. This means that $v_i$ is the $i$-th generalized eigenvector of the matrix pair $(H,\Sigma^{-1})$ and the associated eigenvalue is $\lambda_i=a_i^2$.
 Therefore $P_r$ satisfies $P_r = \Sigma^{1/2} V_rV_r^T \Sigma^{-1/2} = (\sum_{i=1}^r v_iv_i^T ) \Sigma^{-1}$ as in \eqref{eq:OptimalPr} and 
 $\trace( \Sigma(I_d-P_r^T) H (I_d-P_r) ) = \|A-A_r\|_F^2 = \| U(D-D_r)V^T\|^2_F=\lambda_{r+1}+\hdots+\lambda_d$ as in \eqref{eq:OptimalPr_residual}.

\subsection{Proof of Proposition \ref{prop:KL_2}}\label{proof:prop:KL_2}

 The \emph{trace duality} property allows writing $\trace(AB)\leq \|A\|\trace(B)$ for any symmetric positive-semidefinite matrices $A,B\in\mathbb{R}^{d\times d}$, where $\|A\|=\sup\{ |x^TAx| ~,~ x\in\mathbb{R}^d \text{ s.t. } \|x\|_2=1 \}$ denotes the spectral norm of $A$. With the choice $A=H$ and $B=(I_d-Q_r)\Sigma(I_d-Q_r)^T$ we can write
 \begin{align*}
  \trace\big( \Sigma(I_d-Q_r)^T H (I_d-Q_r) \big) 
  &= \trace\big( H (I_d-Q_r) \Sigma(I_d-Q_r)^T \big) \\
  & \leq \|H\| \trace\big( (I_d-Q_r) \Sigma(I_d-Q_r)^T \big) \\
  &= \|H\| \, \mathbb{E}\big( \| (X-m) - Q_r(X-m) \|_2^2 \big) ,
 \end{align*}
 for any projector $Q_r$. Let $Q_r$ be a solution to \eqref{eq:KLmotivation} and $P_r$ be a minimizer of $P_r\mapsto \trace\big( \Sigma(I_d-P_r)^T H (I_d-P_r) \big)$. By Propositions \ref{prop:BoundVVfunction} and \ref{prop:OptimalPr} we can write
 \begin{align*}
  \| f - \mathbb{E}_\mu(f|\sigma(P_r)) \|_{\mathcal{H}}^2 
  \leq \sum_{i=r+1}^d \lambda_i 
  &= \trace\big( \Sigma(I_d-P_r)^T H (I_d-P_r) \big) \\
  &\leq \trace\big( \Sigma(I_d-Q_r)^T H (I_d-Q_r) \big) \\
  &\leq \|H\| \, \mathbb{E}\big( \| (X-m) - Q_r(X-m) \|_2^2 \big)\\
  &= \|H\| \sum_{i=r+1}^d \sigma_i^2 \,.
 \end{align*}
 To conclude the proof, it remains to show that $\|H\|\leq L^2$. 
 Because $f$ is continuously differentiable we can write
 $
  f(x+h)=f(x) + \nabla f(x)h+ \text{o}(\|h\|_2)
 $
 for any $x,h\in\mathbb{R}^d$.
 Also, because $f$ is Lipschitz we have
 \begin{align*}
  \| \nabla f(x) h \|_V &= \| f(x+h) - f(x) + \text{o}(\|h\|_2) \|_V \\
  &\leq \| f(x+h) - f(x)  \|_V + \text{o}(\|h\|_2) \\
  &\leq  L\|h\|_2+ \text{o}(\|h\|_2) \, ,
 \end{align*}
 for any $x,h\in\mathbb{R}^d$. 
 Replacing $h$ by $ty$ where $t>0$ and $\|y\|_2=1$, and dividing by $t$ we obtain 
 $
  \| \nabla f(x) y \|_V \leq L + \text{o}(1) \underset{t\rightarrow 0}{\longrightarrow} L
 $
 for any $\|y\|_2=1$. Thus we have
 \begin{align*}
  \|H \| 
  = \sup_{y\in\mathbb{R}^d,\|y\|_2=1} |y^THy| 
  &= \sup_{y\in\mathbb{R}^d,\|y\|_2=1} \int_{\mathbb{R}^d} \| \nabla f(x) y \|_V^2  \,\mu( \mathrm{d}x) \\
  &\leq  \sup_{y\in\mathbb{R}^d,\|y\|_2=1} L^2\|y\|_2^2 = L^2 \,,
 \end{align*}
 which concludes the proof.

\section*{Acknowledgments}
This material was based upon work partially supported by the National Science Foundation under Grant DMS-1638521 to the Statistical and Applied Mathematical Sciences Institute. Any opinions, findings, and conclusions or recommendations expressed in this material are those of the authors and do not necessarily reflect the views of the National Science Foundation.
O.~Zahm and Y.~Marzouk gratefully acknowledge support from the DARPA EQUiPS program. O.~Zahm, C.~Prieur, and Y.~Marzouk also gratefully acknowledge support from the Inria associate team
UNQUESTIONABLE (UNcertainty QUantification is ESenTIal for OceaNic \& Atmospheric flow proBLEms).

\bibliographystyle{siamplain}
\bibliography{references}

\begin{thebibliography}{10}

\bibitem{Adragni09}
{\sc K.~P. Adragni and R.~D. Cook}, {\em Sufficient dimension reduction and
  prediction in regression}, Philosophical Transactions of the Royal Society of
  London A: Mathematical, Physical and Engineering Sciences, 367 (2009),
  pp.~4385--4405, \url{https://doi.org/10.1098/rsta.2009.0110}.

\bibitem{barreda2007some}
{\sc L.~Barreda, A.~Gannoun, and J.~Saracco}, {\em Some extensions of
  multivariate sliced inverse regression}, Journal of Statistical Computation
  and Simulation, 77 (2007), pp.~1--17.

\bibitem{boucheron2013concentration}
{\sc S.~Boucheron, G.~Lugosi, and P.~Massart}, {\em {Concentration
  inequalities: A nonasymptotic theory of independence}}, Oxford university
  press, 2013.

\bibitem{charrier2013finite}
{\sc J.~Charrier, R.~Scheichl, and A.~L. Teckentrup}, {\em Finite element error
  analysis of elliptic pdes with random coefficients and its application to
  multilevel monte carlo methods}, SIAM Journal on Numerical Analysis, 51
  (2013), pp.~322--352.

\bibitem{Chen1982}
{\sc L.~H. Chen}, {\em An inequality for the multivariate normal distribution},
  Journal of Multivariate Analysis, 12 (1982), pp.~306--315,
  \url{https://doi.org/10.1016/0047-259X(82)90022-7}.

\bibitem{Cohen2012b}
{\sc A.~Cohen, I.~Daubechies, R.~DeVore, G.~Kerkyacharian, and D.~Picard}, {\em
  {Capturing Ridge Functions in High Dimensions from Point Queries}},
  Constructive Approximation, 35 (2012), pp.~225--243.

\bibitem{Constantine2015}
{\sc P.~G. Constantine}, {\em Active Subspaces: Emerging Ideas for Dimension
  Reduction in Parameter Studies}, Society for Industrial and Applied
  Mathematics, Philadelphia, 2015,
  \url{https://doi.org/10.1137/1.9781611973860}.

\bibitem{constantine2017global}
{\sc P.~G. Constantine and P.~Diaz}, {\em {Global sensitivity metrics from
  active subspaces}}, Reliability Engineering \& System Safety, 162 (2017),
  p.~1–13, \url{https://doi.org/10.1016/j.ress.2017.01.013}.

\bibitem{Constantine2017a}
{\sc P.~G. Constantine and A.~Doostan}, {\em Time-dependent global sensitivity
  analysis with active subspaces for a lithium ion battery model}, Statistical
  Analysis and Data Mining: The ASA Data Science Journal, 10 (2017),
  pp.~243--262, \url{https://doi.org/10.1002/sam.11347}.

\bibitem{Constantine2014a}
{\sc P.~G. Constantine, E.~Dow, and Q.~Wang}, {\em Active subspace methods in
  theory and practice: Applications to kriging surfaces}, SIAM Journal on
  Scientific Computing, 36 (2014), pp.~A1500--A1524,
  \url{https://doi.org/10.1137/130916138}.

\bibitem{cookweis1991}
{\sc R.~Cook and S.~Weisberg}, {\em Discussion of 'sliced inverse regression
  for dimension reduction'}, Journal of the American Statistical Association,
  86 (1991), pp.~328--332.

\bibitem{Cook1998}
{\sc R.~D. Cook}, {\em Regression Graphics: Ideas for Studying Regressions
  through Graphics}, John Wiley \& Sons, Inc., New York, 1998,
  \url{https://doi.org/10.1002/9780470316931}.

\bibitem{Diaconis1984}
{\sc P.~Diaconis and M.~Shahshahani}, {\em On nonlinear functions of linear
  combinations}, SIAM Journal on Scientific and Statistical Computing, 5
  (1984), pp.~175--191, \url{https://doi.org/10.1137/0905013}.

\bibitem{Ern2004}
{\sc A.~Ern and J.-L. Guermond}, {\em {Theory and Practice of Finite
  Elements}}, vol.~159 of Applied Mathematical Sciences, Springer New York, New
  York, NY, 2004.

\bibitem{Fornasier2012}
{\sc M.~Fornasier, K.~Schnass, and J.~Vybiral}, {\em {Learning functions of few
  arbitrary linear parameters in high dimensions}}, Foundations of
  Computational Mathematics, 12 (2012), pp.~229--262.

\bibitem{frauenfelder2005finite}
{\sc P.~Frauenfelder, C.~Schwab, and R.~A. Todor}, {\em Finite elements for
  elliptic problems with stochastic coefficients}, Computer methods in applied
  mechanics and engineering, 194 (2005), pp.~205--228.

\bibitem{Friedman1981}
{\sc J.~H. Friedman and W.~Stuetzle}, {\em Projection pursuit regression},
  Journal of the American Statistical Association, 76 (1981), pp.~817--823,
  \url{https://doi.org/10.1080/01621459.1981.10477729}.

\bibitem{gamboa:hal-00800847}
{\sc F.~Gamboa, A.~Janon, T.~Klein, and A.~Lagnoux}, {\em Sensitivity indices
  for multivariate outputs}, Comptes Rendus Mathematique, 351 (2013), pp.~307
  -- 310.

\bibitem{gamboa:hal-00881112}
{\sc F.~Gamboa, A.~Janon, T.~Klein, A.~Lagnoux, et~al.}, {\em Sensitivity
  analysis for multidimensional and functional outputs}, Electronic Journal of
  Statistics, 8 (2014), pp.~575--603.

\bibitem{hartapproximation}
{\sc J.~Hart and P.~Gremaud}, {\em An approximation theoretic perspective of
  sobol'indices with dependent variables}, International Journal for
  Uncertainty Quantification,
  \url{https://doi.org/10.1615/Int.J.UncertaintyQuantification.2018026498}.

\bibitem{Haykin1999}
{\sc S.~Haykin}, {\em Neural Networks: A Comprehensive Foundation}, Prentice
  Hall, Upper Saddle River, NJ, 2nd~ed., 1999.

\bibitem{homma1996importance}
{\sc T.~Homma and A.~Saltelli}, {\em Importance measures in global sensitivity
  analysis of nonlinear models}, Reliability Engineering \& System Safety, 52
  (1996), pp.~1--17.

\bibitem{Huber1985}
{\sc P.~J. Huber}, {\em Projection pursuit}, The Annals of Statistics, 13
  (1985), pp.~435--475, \url{http://www.jstor.org/stable/2241175}.

\bibitem{iooss2015review}
{\sc B.~Iooss and P.~Lema{\^\i}tre}, {\em A review on global sensitivity
  analysis methods}, in Uncertainty management in simulation-optimization of
  complex systems, Springer, 2015, pp.~101--122.

\bibitem{Jefferson2017}
{\sc J.~L. Jefferson, R.~M. Maxwell, and P.~G. Constantine}, {\em Exploring the
  sensitivity of photosynthesis and stomatal resistance parameters in a land
  surface model}, Journal of Hydrometeorology, 18 (2017), pp.~897--915,
  \url{https://doi.org/10.1175/JHM-D-16-0053.1}.

\bibitem{ji2018shared}
{\sc W.~Ji, J.~Wang, O.~Zahm, Y.~M. Marzouk, B.~Yang, Z.~Ren, and C.~K. Law},
  {\em Shared low-dimensional subspaces for propagating kinetic uncertainty to
  multiple outputs}, Combustion and Flame, 190 (2018), pp.~146--157.

\bibitem{Jolliffe2002}
{\sc I.~T. Jolliffe}, {\em Principal Component Analysis}, Springer, New York,
  2nd~ed., 2002, \url{https://doi.org/10.1007/b98835}.

\bibitem{kucherenko2009monte}
{\sc S.~Kucherenko, M.~Rodriguez-Fernandez, C.~Pantelides, and N.~Shah}, {\em
  Monte carlo evaluation of derivative-based global sensitivity measures},
  Reliability Engineering \& System Safety, 94 (2009), pp.~1135--1148.

\bibitem{kucherenko2009derivative}
{\sc S.~Kucherenko and I.~M. Sobol}, {\em Derivative based global sensitivity
  measures and their link with global sensitivity indices}, Mathematics and
  Computers in Simulation, 79 (2009), pp.~3009--3017.

\bibitem{kucherenko2016derivative}
{\sc S.~Kucherenko and S.~Song}, {\em Derivative-based global sensitivity
  measures and their link with sobol sensitivity indices}, in Monte Carlo and
  Quasi-Monte Carlo Methods, Springer, 2016, pp.~455--469.

\bibitem{lamboni2013derivative}
{\sc M.~Lamboni, B.~Iooss, A.-L. Popelin, and F.~Gamboa}, {\em Derivative-based
  global sensitivity measures: general links with sobol’indices and numerical
  tests}, Mathematics and Computers in Simulation, 87 (2013), pp.~45--54.

\bibitem{lamboni2011multivariate}
{\sc M.~Lamboni, H.~Monod, and D.~Makowski}, {\em Multivariate sensitivity
  analysis to measure global contribution of input factors in dynamic models},
  Reliability Engineering \& System Safety, 96 (2011), pp.~450--459.

\bibitem{li1991sliced}
{\sc K.-C. Li}, {\em Sliced inverse regression for dimension reduction},
  Journal of the American Statistical Association, 86 (1991), pp.~316--327.

\bibitem{li1992principal}
{\sc K.-C. Li}, {\em On principal hessian directions for data visualization and
  dimension reduction: Another application of stein's lemma}, Journal of the
  American Statistical Association, 87 (1992), pp.~1025--1039.

\bibitem{li2003dimension}
{\sc K.-C. Li, Y.~Aragon, K.~Shedden, and C.~Thomas~Agnan}, {\em Dimension
  reduction for multivariate response data}, Journal of the American
  Statistical Association, 98 (2003), pp.~99--109.

\bibitem{Lukaczyk2014}
{\sc T.~W. Lukaczyk, P.~Constantine, F.~Palacios, and J.~J. Alonso}, {\em
  Active subspaces for shape optimization}, in 10th AIAA Multidisciplinary
  Design Optimization Conference, 2014,
  \url{https://doi.org/10.2514/6.2014-1171}.

\bibitem{Marzouk2009}
{\sc Y.~M. Marzouk and H.~N. Najm}, {\em Dimensionality reduction and
  polynomial chaos acceleration of bayesian inference in inverse problems},
  Journal of Computational Physics, 228 (2009), p.~1862–1902,
  \url{https://doi.org/10.1016/j.jcp.2008.11.024}.

\bibitem{Mayer2015}
{\sc S.~Mayer, T.~Ullrich, and J.~Vyb{\'{i}}ral}, {\em {Entropy and Sampling
  Numbers of Classes of Ridge Functions}}, vol.~42, Springer US, 2015.

\bibitem{Mikosch1998}
{\sc T.~Mikosch and O.~Kallenberg}, {\em {Foundations of Modern Probability}},
  Journal of the American Statistical Association, 93 (1998), p.~1243,
  \url{https://doi.org/10.2307/2669881}.

\bibitem{nobile2016adaptive}
{\sc F.~Nobile, L.~Tamellini, F.~Tesei, and R.~Tempone}, {\em An adaptive
  sparse grid algorithm for elliptic pdes with lognormal diffusion
  coefficient}, in Sparse Grids and Applications-Stuttgart 2014, Springer,
  2016, pp.~191--220.

\bibitem{nourdin2009second}
{\sc I.~Nourdin, G.~Peccati, and G.~Reinert}, {\em Second order {P}oincar{\'e}
  inequalities and {CLTs} on {Wiener} space}, Journal of Functional Analysis,
  257 (2009), pp.~593--609.

\bibitem{Pinkus2015}
{\sc A.~Pinkus}, {\em Ridge Functions}, Cambridge University Press, Cambridge,
  2015, \url{https://doi.org/10.1017/CBO9781316408124}.

\bibitem{plessix2006review}
{\sc R.-E. Plessix}, {\em {A review of the adjoint-state method for computing
  the gradient of a functional with geophysical applications}}, Geophysical
  Journal International, 167 (2006), pp.~495--503.

\bibitem{roustant2017poincare}
{\sc O.~Roustant, F.~Barthe, B.~Iooss, et~al.}, {\em Poincar{\'e} inequalities
  on intervals--application to sensitivity analysis}, Electronic journal of
  statistics, 11 (2017), pp.~3081--3119.

\bibitem{russi2010uncertainty}
{\sc T.~M. Russi}, {\em Uncertainty quantification with experimental data and
  complex system models}, PhD thesis, UC Berkeley, 2010.

\bibitem{Saltelli2008}
{\sc A.~Saltelli, M.~Ratto, T.~Andres, F.~Campolongo, J.~Cariboni, D.~Gatelli,
  M.~Saisana, and S.~Tarantola}, {\em Global Sensitivity Analysis. The Primer},
  John Wiley \& Sons, Inc., New York, 2008,
  \url{https://doi.org/10.1002/9780470725184}.

\bibitem{saltar04}
{\sc A.~Saltelli, S.~Tarantola, F.~Campolongo, and M.~Ratto}, {\em Sensitivity
  analysis in practice: A guide to assessing scientific models}, Wiley, 2004.

\bibitem{samarov1993exploring}
{\sc A.~M. Samarov}, {\em Exploring regression structure using nonparametric
  functional estimation}, Journal of the American Statistical Association, 88
  (1993), pp.~836--847.

\bibitem{saracco2005asymptotics}
{\sc J.~Saracco}, {\em Asymptotics for pooled marginal slicing estimator based
  on sir$\alpha$ approach}, Journal of multivariate Analysis, 96 (2005),
  pp.~117--135.

\bibitem{Schwab2006}
{\sc C.~Schwab and R.~A. Todor}, {\em {K}arhunen--{L}o\`{e}ve approximation of
  random fields by generalized fast multipole methods}, Journal of
  Computational Physics, 217 (2006), pp.~100--122,
  \url{https://doi.org/10.1016/j.jcp.2006.01.048}.

\bibitem{sobol1993sensitivity}
{\sc I.~M. Sobol}, {\em Sensitivity estimates for nonlinear mathematical
  models}, Mathematical modelling and computational experiments, 1 (1993),
  pp.~407--414.

\bibitem{sobol2001global}
{\sc I.~M. Sobol}, {\em Global sensitivity indices for nonlinear mathematical
  models and their monte carlo estimates}, Mathematics and computers in
  simulation, 55 (2001), pp.~271--280.

\bibitem{zahm2018certified}
{\sc O.~Zahm, T.~Cui, K.~Law, A.~Spantini, and Y.~Marzouk}, {\em Certified
  dimension reduction in nonlinear bayesian inverse problems},
  arXiv:1807.03712,  (2018).

\bibitem{zhu2010dimension}
{\sc L.-P. Zhu, L.-X. Zhu, and S.-Q. Wen}, {\em On dimension reduction in
  regressions with multivariate responses}, Statistica Sinica,  (2010),
  pp.~1291--1307.

\end{thebibliography}
\end{document}